\documentclass[11pt]{article}
\usepackage{amsmath,amsfonts,amscd,bezier}
\usepackage{latexsym,amssymb}
\usepackage[utf8]{inputenc} 
\usepackage[T1]{fontenc}    
\usepackage[english]{babel}
\usepackage[usenames]{color}

\newtheorem{lemma}{Lemma}[section]
\newtheorem{theorem}[lemma]{Theorem}
\newtheorem{remark}[lemma]{Remark}
\newtheorem{proposition}[lemma]{Proposition}
\newtheorem{corollary}[lemma]{Corollary}

\newtheorem{example}[lemma]{Example}

\newcommand{\cqd}{{\hfill $\rule{2mm}{2mm}$}\vspace{1cm}}
\newcommand{\Dem}{\noindent{\sc Proof:\ \ }}

\begin{document}
	
	\title{\bf Dicritical foliations and semiroots of plane branches}
	
	\author{\sc Nuria Corral, Marcelo E. Hernandes and Maria Elenice R. Hernandes \thanks{The authors are supported by the Spanish research projects PID2019-105621GB-I00 and PID2022-139631NB-I00 funded by the Agencia Estatal de Investigaci\'on - Ministerio de Ciencia e Innovaci\'on. The second author is also supported by the research project 303638/2020-6 funded by CNPq-Brazil.}
	}
	
	\date{}

	\maketitle

	\begin{abstract}
		In this work we describe dicritical foliations in $(\mathbb{C}^2,0)$ at a triple point of the resolution dual graph of an analytic plane branch $\mathcal{C}$ using its semiroots. In particular, we obtain a constructive method to present a one-parameter family $\mathcal{C}_{u}$ of separatrices for such foliations. As a by-product we relate the contact order between a special member of $\mathcal{C}_{u}$ and $\mathcal{C}$ with analytic discrete invariants of plane branches.
		
	\end{abstract}
	
	\noindent {2020 AMS Classification: } 32S65, 14H20, 32S10 \\

\section{Introduction}

The aim of this work is to describe a construction of foliations in $(\mathbb{C}^2,0)$ with a dicritical component in one of the bifurcation divisor of the reduction of singularities of an irreducible plane curve (branch).

Previous works deal with the construction of dicritical foliations. For instance, in \cite{Can-C-2011}, it is proved the existence of absolutely dicritical foliations for any configuration of the exceptional divisor, that is, given a morphism $\sigma: M \to (\mathbb{C}^2,0)$ composition of a finite number of punctual blow-ups, there exists a germ of foliation $\mathcal{F}$ in $(\mathbb{C}^2,0)$ such that the transformed foliation $\sigma^*\mathcal{F}$ is completely transversal to the exceptional divisor $\sigma^{-1}(0)$. Moreover, the foliation $\mathcal{F}$ has a meromorphic first integral. In \cite{Can-C-2006}, the authors present a way to construct logarithmic dicritical foliations (weak logarithmic models)  which share some properties with a given foliation. More precisely, a weak logarithmic model $\mathcal{L}$ for a foliation $\mathcal{F}$ is a logarithmic foliation such that the reduction of the singularities of $\mathcal{L}$ is longer than the one of $\mathcal{F}$ and coincides with it outside a “escape set” of non-singular points for $\mathcal{F}$ placed at dicritical components. Note that, this escape set depends on the analytic type of the curves defined by the equations used to write the 1-form which gives the logarithmic foliation as shown in \cite[Example 19]{Can-C-2006}.

Our approach here is different from the works previously mentioned. The main tools of our construction are the concept of the semiroots of a branch $\mathcal{C}$, which codify part of the topological data of the curve, and a result concerning a special way to express a holomorphic 1-form in $\Omega^1$ (Azevedo's Lemma).

Semiroots of a plane branch $\mathcal{C}$ are particular branches that allow us to determine the topological class of the curve. Zariski in \cite{Zariski} considered semiroots in order to relate the characteristic exponents of $\mathcal{C}$ and the minimal generators of the value semigroup associated to the branch. In \cite{Abh}, Abhyankar and Moh introduced particular semiroots (approximate roots) that can be used in an effective criterion of irreducibility of elements in $\mathbb{C}\{x,y\}$.

The other ingredient is the Azevedo's Lemma (see \cite{azevedo}, Chapter 5, Proposition 2): given $n,m\in\mathbb{Z}_{>0}$, any $1$-form $\omega \in \Omega^1$ can be expressed as
$\omega =H_1\cdot (nxdy-mydx)+dH_2,$ with $H_1,H_2\in\mathbb{C}\{x,y\}$. This particular way to express a $1$-form has been used by other authors on topics related to plane curves, vector fields, etc. For instance, Loray in \cite{loray} presents normal forms for cuspidal singularities of analytic vector fields that correspond a particular case of Azevedo's Lemma. Bayer and Hefez (see \cite{bayer}) use such expression as a tool to describe (up to analytic equivalence) plane branches such that the Milnor and Tjurina numbers differ by one or two.


In this work, we consider $\mathcal{C}:=\{F=0\}$ a plane branch, where
$F\in\mathbb{C}\{x\}[y]$ is a Weierstrass polynomial, with value semigroup $\Gamma$ minimally generated by $\{ v_0,\ldots ,v_g\}$. A set $\{F_0,F_1,\ldots ,F_{g+1}:=F\}$ is an extended system of
semiroots of $F$, if $F_i\in\mathbb{C}\{x\}[y]$ is monic, $deg_yF_i=\frac{v_0}{GCD(v_0,\ldots ,v_i)}$ for $1 \leq i \leq g$ and  $\dim_{\mathbb{C}}\frac{\mathbb{C}\{x,y\}}{\langle F,F_i\rangle}=v_i$ for $0\leq i\leq g$ (see Section 3).

For each pair $(i,j)$ with $0 \leq i < j\leq g$ we set $\omega_{ij}=v_{i}F_idF_j - v_{j}F_j dF_i\in\Omega^1$ and consider the singular foliation $\mathcal{F}_{\omega}$ defined by
\begin{equation}
\label{Introd-Omega}
\omega= H_1 \omega_{ij} + dH_2,
\end{equation} where $H_1,H_2 \in \mathbb{C}\{x,y\}\setminus \langle F\rangle$. Note that, for $i=0$, $j=1$, $v_0=n$ and $v_1=m$ we get the expression given in Azevedo's Lemma.

In Section 3, we present the main results of this paper, Theorem \ref{TeoSepBq} and Corollary \ref{CorSep}. We give a necessary and sufficient condition to assure that the foliation $\mathcal{F}_{\omega}$, with $\omega$ as in (\ref{Introd-Omega}), has a dicritical component in the last triple point of the resolution dual graph of $\mathcal{C}$. This condition is given in terms of  the intersection multiplicities of $F$ with $H_1$ and $H_2$.
Moreover, in the proof of this theorem we have a constructive method to describe a family of parameterizations for the separatrices in such dicritical component up to the desired order as illustrated in Example \ref{exemplo-sep}. The presented method in Theorem \ref{TeoSepBq} does not make use of blowing up which is normally considered to present dicritical foliations.

Dicritical foliations and analytic invariants of irreducible plane curves are closely related. More precisely, given a parametrization $\varphi(t)=(x(t),y(t))$ of the irreducible plane curve $\mathcal{C}$ and a 1-form $\omega$, we define the value of $\omega$ by $\nu(\omega)=ord_t(\varphi^*(\omega))+1$ where $\varphi^*(\omega)=\omega(\varphi(t))$. The set of differential values $\Lambda$ of $\mathcal{C}$ is given by
 $$\Lambda =\{ \nu(\omega) \ : \ \omega \in \Omega^1\}.$$ This analytic invariant is one of the main ingredients in the analytic classification of branches (see \cite{HH,handbook}). When $\omega=0$ defines a foliation $\mathcal{F}=\mathcal{F}_\omega$  in $(\mathbb{C}^2,0)$ and $\mathcal{C}$ is not an invariant curve (separatrix) of $\mathcal{F}$, the value $\nu(\omega)-1$ coincides with the tangency order $\tau_0(\mathcal{F},\mathcal{C})$ of the foliation $\mathcal{F}$  with the curve $\mathcal{C}$ (see \cite{Cam-LN-S,Can-C-M-2019}). If we consider a hamiltonian 1-form $\omega=dg$, with $g \in \mathbb{C}\{x,y\}$ a non unit, then $\nu(dg)=(\mathcal{C},\mathcal{D})_0$, where $(\mathcal{C},\mathcal{D})_0$ denotes the intersection multiplicity at the origin of the curves $\mathcal{C}$ and $\mathcal{D}$, with $\mathcal{D}:=\{g=0\}$. Hence, we have that $\Gamma \setminus \{0\} \subseteq \Lambda$ where $\Gamma$ is the value semigroup  associated to the curve $\mathcal{C}$.
 Moreover, there exists a finite subset $L=\{\ell_1,\ldots, \ell_k\} \subset \Lambda$ such that any $\ell \in \Lambda$ can be expressed as $\ell=\ell_i + \gamma$ for some $\gamma \in \Gamma$ and $\ell_i \in L$, that is, the set $\Lambda$ is a finitely generated $\Gamma$-monomodule.

 Let us consider a set of 1-forms $\{\omega_1,\ldots, \omega_k\}$ such that $\nu(\omega_i)=\ell_i$. If the curve $\mathcal{C}$ has only one Puiseux pair, then the foliations defined by the 1-forms $\omega_i=0$ are dicritical in the triple point of the resolution dual graph of the curve $\mathcal{C}$ (see \cite{nuria} where properties of these 1-forms are described).

From the results in \cite{Can-C-M-2019,Cor-2024}, we have that, if the foliation $\mathcal{F}$ defined by $\omega=0$ is a non-dicritical second type foliation (see \cite{Mat-S}), then $\tau_0(\mathcal{F},\mathcal{C})=\nu(\omega)-1=(S_\mathcal{F},\mathcal{C})_0-1$, where $S_\mathcal{F}$ is the curve of separatrices of $\mathcal{F}$. Thus, if $\omega$ is a 1-form such that $\nu(\omega) \in \Lambda \setminus \Gamma$,  the foliation defined by $\omega=0$ is either dicritical or it is not a second type foliation.

In Section 4, we explore $1$-forms expressed as $\omega= H_1 \omega_{ij} + dH_2$ and their connection with the analytic invariant $\Lambda$. We show that the value of $\omega$ is related with the contact between the branch and a special separatrix $\mathcal{C}_{{\star}}$ of $\mathcal{F}_{\omega}$ (see Theorem \ref{cont-value}). In particular, for curves $\mathcal{C}$ with semigroup $\langle v_0,v_1\rangle$, we show that the set $\Lambda$ can be determined using dicritical foliations defined by $H_1 \omega_{01} + dH_2$, or equivalently, by $H_1$ and the special separatrix $\mathcal{C}_{{\star}}$ (see Corollary \ref{Lambda-g-1}). The separatrix $\mathcal{C}_{{\star}}$ is closely related to the concept of {\it analytic semiroot} introduced by Cano, Corral and Senovilla-Sanz in \cite{nuria} were, as we mentioned before, geometrical properties are presented for $\Lambda$. In addition, Proposition \ref{oziel-generalized} shows how to compute the Zariski invariant $\lambda$ of $\mathcal{C}$, that is $\lambda=\min(\Lambda\setminus\Gamma)-v_0$, considering dicritical foliations in the first triple point of the resolution dual graph of $\mathcal{C}$ that extends a result by G\'omez-Mart\'inez presented for branches with value semigroup minimally generated by two elements.

\vspace{0.5cm}

{\bf Acknowledgements:} We are grateful to Prof. Felipe Cano for all conversations and suggestions.
The last two authors would like to express their gratitude to the ECSING-AFA group from Universidad de Valladolid for its hospitality during the year 2013, when the seminal ideas of this work emerged.

\section{Notations}

In this section we present some classic notations. For the results about Plane Curve Theory and Foliation Theory we indicate \cite{hefez} and \cite{cano}, respectively. We denote by $\mathbb{C}\{x,y\}$ the absolutely convergent power series ring at the origin in $\mathbb{C}^2$.

A germ of an analytic plane curve $\mathcal{C}_F$ in $(\mathbb{C}^2,0)$
is the (germ of) zero set of a reduced element $F \in
\mathbb{C}\{x,y\}$ in a neighborhood at the origin. Without loss of generality (by a change of coordinates) we can
consider $F\in\mathbb{C}\{x\}[y]$ a Weierstrass polynomial
$F(x,y)=y^n+\sum_{i=1}^{n}A_i(x)y^{n-i}$ where $n$ is the multiplicity of $F$, denoted by $mult(F)$.

If $F$
is irreducible, we can assume that $y=0$ is the tangent cone of the branch
$\mathcal{C}_F$ and this implies that $mult(A_i(x))>i$ for $1\leq i\leq n$. By
Newton-Puiseux theorem we can obtain
$\eta\left (x^{\frac{1}{n}}\right )=\sum_{k>n}c_kx^{\frac{k}{n}}\in\mathbb{C}\left \{x^{\frac{1}{n}}\right \}$
such that $F\left (x,\eta\left (x^{\frac{1}{n}}\right )\right )=0$ and the set of roots of
$F$ (in a neighborhood at the origin) is
$\left \{\eta\left (\alpha\cdot x^{\frac{1}{n}}\right );\ \alpha\in U_n\right \}$, where
$U_n=\{\alpha\in\mathbb{C};\ \alpha^n=1\}$. In particular, we have
\begin{equation}
\label{DecRootUnity} F(x,y)=\prod_{\alpha\in
	U_n}\left (y-\eta\left (\alpha\cdot x^{\frac{1}{n}}\right )\right ).
\end{equation}

By a Tschirnhausen transformation, i.e. by the change of coordinates $(x,y)\rightarrow \left (x,y-\frac{1}{n}A_1(x)\right )$, we can assume that $A_1(x)=0$, or equivalently $c_k=0$ for all $k\equiv 0\mod n$ in $\eta\left (x^{\frac{1}{n}}\right )$.

Putting $t=x^{\frac{1}{n}}$ we obtain a Puiseux parameterization for
$\mathcal{C}_F$:
\begin{equation}\label{PuiseuxParam}
\varphi(t)=\left (t^{n},\sum_{k\geq\beta_1}c_kt^k\right ),
\end{equation} where $\beta_1=\min\{k;\ k\not\equiv 0\mod n\ \mbox{and}\ c_k\neq 0\}$. Moreover, we will
assume that such parameterization is primitive, that is, $\varphi(t)$
can not be reparameterized by a power of a new variable or equivalently
the greatest common divisor of all exponents in $\varphi(t)$ is
equal to $1$.

In what follows we consider plane branches, that is, plane curves defined by an irreducible Weierstrass polynomial as (\ref{DecRootUnity}).

There are two sequences $(e_i)$ and
$(\beta_i)$ of integers associated to $\mathcal{C}_F$ and obtained
by any Puiseux parameterization of $\mathcal{C}_F$:
$$\begin{array}{l}
\beta_0=e_0=n;\vspace{0.15cm} \\
\beta_j=min\{i; \ i \not\equiv 0 \ mod \ e_{j-1} \ \mbox{and}\ c_i \neq 0 \}; \vspace{0.15cm} \\
e_{j}=GCD(e_{j-1},\beta_j)=GCD(\beta_0, \ldots, \beta_j).
\end{array}
$$
The elements in the increasing finite sequence $(\beta_i)_{i=0}^g$ are called characteristic exponents
associated to the branch and such sequence completely characterizes the topological
type of the curve as an immersed germ in $(\mathbb{C}^2,0)$. The local topology of plane branches can also be determined by the value semigroup
$\Gamma_F$ associated to the curve $\mathcal{C}_F$. More explicitly,
$$\Gamma_F=\{I(F,G);\ G \in \mathbb{C}\{x,y\}\}\subset\overline{\mathbb{N}}:=\mathbb{N}\cup\{\infty\},$$ where $I(F,G)=(\mathcal{C}_F,\mathcal{C}_G)_0$ is
the intersection multiplicity of $\mathcal{C}_F$ and $\mathcal{C}_G$ at the origin
that can be computed by
$$I(F,G)=dim_{\mathbb{C}}\frac{\mathbb{C}\{x,y\}}{\langle F,G\rangle}=ord_t(\varphi^*(G)),$$
and $\varphi^*(G):=G(\varphi(t))$ for a parameterization $\varphi(t)$ of $\mathcal{C}_F$ as (\ref{PuiseuxParam}).

Notice that given $F\in\mathbb{C}\{x\}[y]$ with $deg_y(F)=mult(F)=v_0>1$ any $G\in\mathbb{C}\{x,y\}$ can be expressed, by  Weierstrass Division Theorem, as $G=Q F+H$ with $H\in\mathbb{C}\{x\}[y]$ and $deg_y(H)<v_0$. As $I(F,QF+H)=I(F,H)$ we get
$$\Gamma_{F}=\{I(F,H);\ H\in\mathbb{C}\{x\}[y]\ \mbox{with}\  deg_y(H)<v_0\}.$$

Zariski (in
\cite{Zariski}) showed that the value semigroup $\Gamma_F$  is minimally
generated by the set of integers $\{v_0, v_1, \dots, v_g \}$,
inductively defined  by
\begin{equation}\label{relations0}
v_0=\beta_0=n,\ \ \ v_1=\beta_1\ \ \ \mbox{and}\ \ \ v_i=n_{i-1}v_{i-1}+\beta_i-\beta_{i-1}
\end{equation} or
\begin{equation}\label{relations}
v_i=\sum_{j=0}^{i-2}\frac{e_j-e_{j+1}}{e_{i-1}}\beta_{j+1}+\beta_i
\end{equation} for $i=2, \dots,g$ where $n_0=1$ and $n_i=\frac{e_{i-1}}{e_i}$. It follows from the definition of $n_i$ that $n=n_0\cdot n_1\cdot
\ldots \cdot n_g$. We denote $\Gamma_F=\langle v_0, v_1, \dots, v_g \rangle$ and sometimes it would be convenient to consider $\beta_{g+1}=v_{g+1}=\infty$.

The value semigroup $\Gamma_{F}$ admits a conductor $\mu_F$, that is, $\mu_F+\mathbb{N}\subseteq\Gamma_{F}$ and $\mu_F-1\not\in\Gamma_{F}$. For plane branches, $\mu_F$ coincides with the Milnor number of $\mathcal{C}_F$ and
\begin{equation}\label{milnor}
	\mu_F=\dim_{\mathbb{C}}\frac{\mathbb{C}\{x,y\}}{\langle F_x, F_y\rangle}=\sum_{l=1}^{g}(n_l-1)v_l-v_0+1.
\end{equation}

In this paper we consider germs of holomorphic singular foliations  of codimension one in $(\mathbb{C}^2,0)$ locally given by $\omega=0$,
where
$$\omega=A(x,y)dx+B(x,y)dy\in\Omega^1:=\Omega^1_{\mathbb{C}^2,0}=\mathbb{C}\{x,y\}dx+\mathbb{C}\{x,y\}dy$$
with $A,B \in \mathbb{C}\{x,y\}$, $A(0,0)=B(0,0)=0$ and
$GCD(A,B)=1$. A such foliation will be denoted by $\mathcal{F}_{\omega}$ and its singular locus $Sing(\mathcal{F}_{\omega})$ is
locally given by the common zeros of $A$ and $B$.

An analytic plane branch $\mathcal{C}_F$ defined by $F=0$  is called a separatrix
(or an invariant curve) of a foliation $\mathcal{F}_{\omega}$ if
$\omega \wedge dF=F \cdot G\cdot dx\wedge dy$, where $G\in\mathbb{C}\{x,y\}$. In particular $\mathcal{C}_F \setminus
Sing(\mathcal{F}_{\omega})$ is a leaf of $\mathcal{F}_{\omega}$.

If $\varphi(t)$ is a parameterization of $\mathcal{C}_F$ we can define the $\mathbb{C}$-linear map
\begin{equation}\label{phi-omega}
\begin{array}{cccl}
\varphi^* : & \Omega^1 & \rightarrow & \mathbb{C}\{t\} \\
& \omega=Adx+Bdy & \mapsto & \varphi^*(A)x'(t)+\varphi^*(B)y'(t)
\end{array}
\end{equation} and we have that $\varphi^*(\omega)=0$ if and only if $\frac{\omega\wedge dF}{dx\wedge dy}\in\langle F\rangle$, that is, $\mathcal{C}_F$ is a separatrix of $\mathcal{F}_{\omega}$.

A singular foliation $\mathcal{F}_{\omega}$ is called dicritical if there
is a finite sequence of blowing-ups with nonsingular invariant
centers, such that this process leads to an irreducible component $E_i$ of the exceptional divisor $E$ that
is generically transversal to the strict transform of $\mathcal{F}_{\omega}$.
For codimension one foliations in $(\mathbb{C}^2,0)$, the dicritical
condition is equivalent to the property of having infinitely many
transversal invariant curves through almost any point in $E_i$, or in other words there are infinitely germs of analytic
curves (separatrices) containing the origin and invariant by the
foliation. In this case we say that $\mathcal{F}_{\omega}$ is dicritical in $E_i$ or in the point $Q_i$ of the resolution dual graph $G(\mathcal{C}_F)$ corresponding to $E_i$.

In the next section, we will consider particular plane branches $\mathcal{C}_{F_i}$ such that $I(F,F_i)=v_i$ in order to define dicritical foliations in specific components of the exceptional divisor obtained by the canonical resolution of $\mathcal{C}_F$.

\section{Semiroots and dicritical foliations}

Azevedo, in his thesis (see \cite{azevedo}, Chapter 5, Proposition 2), exhibits a particular way to express any $1$-form in $\Omega^1$ as we present below:

\begin{lemma}[Azevedo]\label{azevedo} Given any $n,m\in\mathbb{Z}_{>0}$ and $\omega\in\Omega^1$, there exist $H_1,H_2\in\mathbb{C}\{x,y\}$ such that
	\begin{equation}\label{w01}
		\omega =H_1\cdot (nxdy-mydx)+dH_2.\end{equation}
\end{lemma}
\Dem The proof is constructive and allow us to obtain $H_1$ and $H_2$ satisfying (\ref{w01}). Given any $\omega=Adx+Bdy \in \Omega^1$ it is possible to prove that there exist $H_1=\sum_{i,j\geq 0}a_{ij}x^iy^j$ and $H_2=\sum_{i,j\geq 0}b_{ij}x^iy^j$ such that $(H_2)_x=A+myH_1$ and $(H_2)_y=B-nxH_1$. To do this, we integrate the first equation in $x$ and substitute in the second one. Thus we obtain a recursive expression to determine the coefficients $a_{ij}$ and $b_{ij}$, and consequently $H_1$ and $H_2$ (see \cite{azevedo}, Chapter 5, Lemma 1 or \cite{bayer}, Proposition 2).
\cqd

In some cases, the expression (\ref{w01}) can be useful to determine separatrices of $\mathcal{F}_{\omega}$ directly from numerical data of $H_1$ and $H_2$. The following example illustrates a such situation.
	
\begin{example}\label{monomial}
 In (\ref{w01}), let us consider $H_1=y^a$ and $H_2=e\cdot x^b$ with $a, b\in\mathbb{Z}_{\geq 0}$, $b\neq 0$, $e\in\mathbb{C}^*$ and $n(b-1)\neq m(a+1)$, that is, $\omega=y^{a}\cdot (nxdy-mydx)+d(e\cdot x^{b})$. It is immediate that $x=0$ is a separatrix of $\mathcal{F}_{\omega}$. Moreover, by some computations, we get that a monomial germ $\varphi(t)=(t^{\alpha},ct^{\beta})$ with $c\neq 0$ parameterizes a separatrix of $\mathcal{F}_{\omega}$ if and only if	$$\alpha=a+1,\ \ \beta=b-1,\ \ \mbox{and}\ \ c^{a+1}=-\frac{(a+1)be}{n(b-1)-m(a+1)}.$$
 Notice that $\varphi(t)$ is not necessarily a primitive parameterization.

	Similarly we can obtain the description of monomial separatrices for foliations defined by $\mathcal{F}_{\omega}$ considering $H_1$ and $H_2$ given by other possible monomials.
\end{example}

In the sequel we will take $1$-forms given by a similar expression in Azevedo's Lemma but considering semiroots of an irreducible Weierstrass polynomial $F\in\mathbb{C}\{x\}[y]$.

Let $\mathcal{C}_F$ be an irreducible plane curve defined by
$F\in\mathbb{C}\{x\}[y]$ with semigroup $\Gamma_F=\langle v_0,\ldots ,v_g\rangle$.
By the minimality of the  generators set $\{v_0,\ldots ,v_g\}$, any element $G\in\mathbb{C}\{x,y\}$ such that $I(F,G)=v_i$
is irreducible. As $y=0$ is the
tangent cone of $\mathcal{C}_F$, it follows that $I(F,x)=v_0$.

A set $\{F_i;\ 1\leq i\leq g+1\}\subset\mathbb{C}\{x\}[y]$ of monic
polynomials satisfying
\begin{enumerate}
	\item[i)] $deg_yF_i=n_0\cdot\ldots\cdot n_{i-1}=\frac{v_0}{e_{i-1}};$
	\item[ii)] $I(F,F_i)=v_{i}$
\end{enumerate}
is called a {\bf system of semiroots} of $F$. We say that
$\{F_0:=x,F_1,\ldots ,F_{g+1}:=F\}$ is an {\bf extended system of
	semiroots} of $F$ and $F_i$ is an {\bf $i$th semiroot}\footnote{Some authors (see \cite{Popescu} for instance) consider the $i$th semiroot of $F$ for $0\leq i\leq g$ as a monic polynomial $F_i\in\mathbb{C}\{x\}[y]$ satisfying $deg_yF_i=\frac{v_0}{e_i}$ and $I(F,F_i)=v_{i+1}$.} of $F$, for $0\leq i\leq g+1$. Moreover, for $i\neq 0$ we have that  the semigroup and
the characteristic exponents of
$\mathcal{C}_{F_i}$ (see \cite{Popescu}) are
\begin{equation}\label{semiroot-data}
\Gamma_{F_i}=\left\langle
\frac{v_0}{e_{i-1}},\ldots ,\frac{v_{i-1}}{e_{i-1}}\right\rangle \hspace{0.7cm} \mbox{and} \hspace{0.7cm} \left \{\frac{\beta_0}{e_{i-1}}, \ \ldots \ ,\ \frac{\beta_{i-1}}{e_{i-1}}\right \}.
\end{equation}

If $\{F_0=x,F_1,\ldots ,F_g,F_{g+1}=F\}$ is an extended
system of semiroots of $F$, then we have that  $\{F_0=x,F_1,\ldots ,F_k,F_{k+1}\}$ is an extended
system of semiroots of $F_{k+1}$ (see \cite{Popescu}).

We can obtain a system of semiroots of $F\in\mathbb{C}\{x\}[y]$ by several ways, for instance
considering the approximate roots introduced by Abhyankar and Moh
(see \cite{Abh} or \cite{Popescu}) or taking representatives for elements in a
minimal Standard Basis of $\frac{\mathbb{C}\{x,y\}}{\langle
	F\rangle}$ (see \cite{basestandard}).

In what follows we will
consider a particular system of semiroots following Zariski's
approach (see \cite{Zariski}) obtained by a parameterization
$\varphi(t)=(t^{\beta_0},\sum_{k \geq \beta_1}c_kt^k)$ of $\mathcal{C}_F$.

Let us denote
$$
	\varphi_{i}(t)=\left (t^{\frac{\beta_0}{e_{i-1}}},\eta_i(t)\right ):=\left(t^{\frac{\beta_0}{e_{i-1}}},\sum_{\beta_1 \leq k <
		\beta_i}c_kt^{\frac{k}{e_{i-1}}}\right),
$$
for $i=1,\ldots, g+1$, where $\varphi_1(t)=(t,0)$.

\begin{proposition}[Zariski, \cite{Zariski}]
	\label{Prop-MinimalPolyn}
		If $F_i\in\mathbb{C}\{x\}[y]$ is the minimal polynomial of $\eta_i(x^{\frac{e_{i-1}}{\beta_0}})$ over $\mathbb{C}((x))$ where $1\leq i\leq g+1$, then $\{F_0:=x,F_1,\ldots ,F_{g+1}:=F\}$ is an extended system of semiroots of $F$. In particular, $\varphi_i$ is a Puiseux parameterization of $\mathcal{C}_{F_i}$, for $i=1, \ldots, g+1$.
\end{proposition}
\Dem Denoting by $m_i=\frac{\beta_0}{e_{i-1}}=\frac{v_0}{e_{i-1}}$ we have that the minimal polynomial
$F_i\in\mathbb{C}\{x\}[y]$  of
$\eta_i\left (x^{\frac{1}{m_i}}\right )$ over $\mathbb{C}((x))$ is given (as in (\ref{DecRootUnity})) by
\begin{equation}\label{Fi}F_i(x,y)=\prod_{\alpha\in
		U_{m_i}}(y-\eta_i(\alpha\cdot x^{\frac{1}{m_i}}))\end{equation} with
$U_{m_i}=\{\alpha\in\mathbb{C};
\alpha^{m_i}=1\}$, for $i=1,\ldots ,g+1$ (see, for instance, \cite{Zariski}). In particular $deg_yF_i=\frac{v_0}{e_{i-1}}$ and $I(F,F_i)=v_i$ (see the proof of Lemma \ref{coef-control}). Therefore, $F_i$ is an $i$th semiroot of $F$.\cqd

In this work we take the extended system of semiroots of
$F$ obtained as in the above proposition and we call it the {\bf canonical system of semiroots} of $F$.

As an immediate consequence of the classical Euclidian division algorithm, it is possible to obtain a decomposition of any element $H \in \mathbb{C}\{x\}[y]$ in terms of a system of semiroots of $F$.

\begin{proposition}[see \cite{Abh} or \cite{Popescu}]
	\label{Prop-DecompH}
	If $\{F_0=x,F_1,\ldots ,F_{g+1}=F\}$ is an extended system of semiroots of $F \in \mathbb{C}\{x\}[y]$ then any $H\in \mathbb{C}\{x\}[y]$ has a unique expansion given by
	\begin{equation}\label{semiroots}
		H=\sum_{\delta=(\delta_0,\ldots ,\delta_{g+1})}u_{\delta}F_0^{\delta_0}F_1^{\delta_1}F_2^{\delta_2}\ldots F_{g+1}^{\delta_{g+1}},
	\end{equation}
	where $u_{\delta}\in\mathbb{C}$,
	\begin{equation}\label{conditions}
		0\leq \delta_i <n_i=\frac{e_{i-1}}{e_i}\ \mbox{for}\ i\in\{1,\ldots ,g\}, 0 \leq \delta_{g+1} \leq \left[ \frac{deg_y H}{deg_y F} \right]
	\end{equation}
and $[r]$ denotes the integral part of $r \in \mathbb{R}$. Moreover, the order in $t$ of the terms
	$$\varphi^*(F_0)^{\delta_0}\cdot\varphi^*(F_1)^{\delta_1}\cdot\varphi^*(F_2)^{\delta_2}\cdot\ldots \cdot\varphi^*(F_{g})^{\delta_{g}}$$ are two by two distinct, where $\varphi$ is a parameterization of $\mathcal{C}_F$.
\end{proposition}

By the previous result, if $H=\sum_{\delta}u_{\delta}F_0^{\delta_0}F_1^{\delta_1}\cdot\ldots\cdot F_{g+1}^{\delta_{g+1}}$ then
$$I(F,H)=\underset{\delta}{\min}\left\{\displaystyle\sum_{i=0}^{g}\delta_iv_i;\ \delta=(\delta_0,\ldots ,\delta_g,0)\ \mbox{with}\ u_{\delta}\neq 0\right\}.$$

\begin{remark}
Notice that the expansion (\ref{semiroots}) is not necessarily a finite sum and it is a bit different of the expansion presented in \cite{Popescu}. In fact, according to Corollary 5.4 of \cite{Popescu} any $H\in \mathbb{C}\{x\}[y]$ has a unique expansion given by a {\em finite} sum
$$
	H=\sum_{(\delta_1,\ldots ,\delta_{g+1})}h_{\delta_1,\ldots,\delta_{g+1}}F_1^{\delta_1}F_2^{\delta_2}\ldots F_{g+1}^{\delta_{g+1}},
$$
with $h_{\delta_1,\ldots,\delta_{g+1}}\in\mathbb{C}\{x\}$ and $\delta_i,\ 1\leq i\leq g+1$, satisfying the conditions (\ref{conditions}). Writing $h_{\delta_1,\ldots ,\delta_{g+1}}=\sum_{\delta=(\delta_0,\ldots ,\delta_{g+1})}u_{\delta}F_0^{\delta_0}$ with $u_{\delta}\in\mathbb{C}$ we get the expansion (\ref{semiroots}). Recall that an $i$th semiroot for us corresponds to an $(i-1)$th semiroot in \cite{Popescu} for $1\leq i\leq g+1$.
\end{remark}

Considering the canonical embedded resolution $\pi: M\rightarrow (\mathbb{C}^2,0)$
of $\mathcal{C}_F$ and $G(\mathcal{C}_F)$ the dual graph associated
to it, we have that the semiroot $F_i$ is a {\it curvette}\footnote{A {\it
curvette} with respect to the component $E_i$ of the exceptional
divisor $E$ is the image in $(\mathbb{C}^2,0)$ of a smooth curve in $M$
meeting $E_i$ transversely in a single point which lies on no other
component of $E$.} with respect to a component of the exceptional
divisor $E$ corresponding to the $i$th-endpoint of
$G(\mathcal{C}_F)$ (see {\sc Figure 1}). In particular, the extended system
of semiroots appears as coordinates in the embedded resolution
process of $\mathcal{C}_F$ (see \cite{Popescu}). We denote by $T_i$ the $i$th triple point in
the dual graph $G(\mathcal{C}_F)$ that appears in the canonical
resolution process, or equivalently, the first triple point after that $F_i$ is desingularized, which we
indicate by $\widetilde{F_i}$.

\begin{center}
	\setlength{\unitlength}{1cm}
	\begin{picture}(10,5)
		\put(0,4){\line(1,0){8}} \put(0,4){\vector(1,1){0.5}}\put(0.6,4.5){$\widetilde{F_0}$}
		\put(0.9,3.9){$\star$}\put(1.1,4.1){$T_1$}
		\put(1,4){\line(0,-1){2}}\put(1,2){\vector(1,1){0.5}}\put(1.6,2.5){$\widetilde{F_1}$}
		\put(2.9,3.9){$\star$}\put(3.1,4.1){$T_2$}
		\put(3,4){\line(0,-1){3}}\put(3,1){\vector(1,1){0.5}}\put(3.6,1.5){$\widetilde{F_2}$}
		\put(3.9,3.9){$\star$}\put(4.1,4.1){$T_3$}
		\put(4,4){\line(0,-1){2}}\put(4,2){\vector(1,1){0.5}}\put(4.6,2.5){$\widetilde{F_3}$}
		\put(5,3.5){$\ldots$}
		\put(5.9,3.9){$\star$}\put(5.6,4.2){$T_{g-1}$}
		\put(6,4){\line(0,-1){1}}\put(6,3){\vector(1,1){0.5}}\put(6.6,3.5){$\widetilde{F}_{g-1}$}
		\put(7.9,3.9){$\star$}\put(7.6,4.2){$T_g$}
		\put(8,4){\line(0,-1){3}}\put(8,1){\vector(1,1){0.5}}\put(8.6,1.5){$\widetilde{F_g}$}
		\put(8,4){\vector(1,1){0.5}}\put(8.6,4.5){$\widetilde{F}_{g+1}$}
	\end{picture}
	
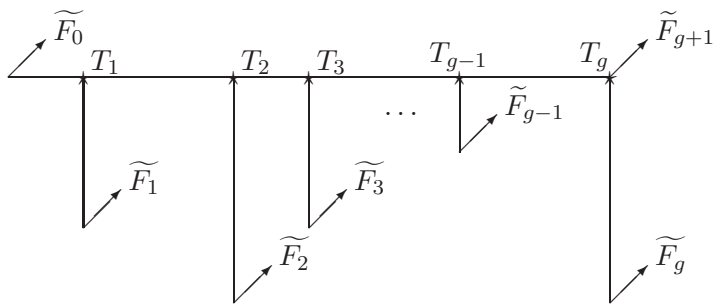
	{\sc Figure 1.} Dual graph for $F_0\cdot F_1\cdot\ldots \cdot F_g\cdot F_{g+1}$.
\end{center}

Given the canonical system of semiroots $\{F_0=x,F_1,\ldots ,F_{g+1}=F\}$ for each $0 \leq i < j\leq g$ we consider $\mathcal{F}_{ij}$ the singular foliation defined by
\begin{equation}\label{wij}\omega_{ij}=v_{i}F_idF_j - v_{j}F_j dF_i.\end{equation}

Notice that $\mathcal{F}_{ij}$ defines the same foliation that $d\left(\frac{F_j^{\alpha_{i}}}{F_i^{\alpha_{j}}}\right)$ where $\alpha_i=\frac{v_i}{GCD(v_i,v_j)}$, $\alpha_j=\frac{v_j}{GCD(v_i,v_j)}$ and therefore $\mathcal{F}_{ij}$ is dicritical with separatrices given by $aF_j^{\alpha_{i}}-b F_i^{\alpha_{j}}=0$ for all $(a:b)\in \mathbb{P}^1_{\mathbb{C}}$.

\begin{example}\label{ex3}
	Let $\mathcal{C}_F$ be the plane branch
	defined by
	$$F=(y^2-x^3-2x^2y+x^4)^3-48x^8(y^2-x^3-2x^2y+x^4)-64x^{11}-64x^{13}$$ with Puiseux parameterization $\varphi (t)=(t^6,t^9+t^{12}+2t^{13})$.
		
	The value semigroup of $\mathcal{C}_F$ is $\Gamma=\langle
	6,9,22\rangle$ and the canonical system of semiroots for it is
		$$F_0=x,\ F_1=y, \ F_2=y^2-x^3-2x^2y+x^4 \ \ \mbox{and}\ \ F_3=F.$$
 We have that
	$\omega_{01}=6xdy-9ydx$ admits separatrices given by
	$ay^2-bx^3=0$, $\omega_{02}=6xdF_2-22F_2dx$ with separatrices $aF_2^3-bx^{11}=0$ and $\omega_{12}=9ydF_2-22F_2dy$ with
	separatrices $aF_2^9-by^{22}=0$, for
	$(a:b)\in\mathbb{P}^{1}_{\mathbb{C}}$.
\end{example}

In what follows we consider $1$-forms given in a particular expression that generalizes (\ref{w01}). More specifically, we take
$\omega=H_1\cdot\omega_{ij}+dH_2$ where $\omega_{ij}$ is given as (\ref{wij}) and admitting that $\omega$ defines a foliation $\mathcal{F}_{\omega}$. We present a simple criterion which assures that $\mathcal{F}_{\omega}$ is dicritical at the $i$th triple point $T_i$ of the dual graph of $\mathcal{C}_{F}$ and we describe a family of separatrices for it.

Our strategy is to consider initially the case $1\leq i<j=g$. The other situations are particular cases of this result by changing $F$ by a semiroot $F_{j+1}$  with $0\leq j< g$.

In order to obtain the results we use some technical lemmas that are presented in Section \ref{technical}.

Given a plane branch $\mathcal{C}_F$ with $F\in\mathbb{C}\{x\}[y]$ as (\ref{DecRootUnity}) and Puiseux parameterization given by $\varphi(t)=\left (t^{\beta_0},\sum_{l\geq\beta_1}c_lt^l\right )$, we consider the family of plane branches $\mathcal{C}_{F_a}$ determined by parameterizations
\begin{equation}\label{family}\varphi_a(t)=\left (t^{\beta_0},\sum_{\beta_1\leq l<\beta_g}c_lt^l+\sum_{l \geq \beta_g}a_lt^l\right )\end{equation}
where $a_l$ are parameters that
can be assume values in $\mathbb{C}$ and $a_{\beta_g}$
nonvanishing. Notice that the coefficient $c_l$ of $t^l$ with
$\beta_1\leq l<\beta_g$ in $\varphi_a(t)$ is precisely the coefficient of
$t^{l}$ in the parameterization $\varphi(t)$ of $\mathcal{C}_F$ and they are considered constant.

It is immediate that $\Gamma_{F_a}=\Gamma_{F}$ and
$\{F_0,F_1,\ldots ,F_{g},F_a\}$ is an extended system of semiroots
for $\mathcal{C}_{F_a}$. Moreover, if $deg_y(H) < deg_y(F)=deg_y(F_a)$ we can write $H=\sum_{\delta}u_{\delta}F_0^{\delta_0}\cdot\ldots\cdot F_{g}^{\delta_g}\in\mathbb{C}\{x\}[y]$ as (\ref{semiroots}) and $I(F,H)=I(F_a,H)$.

\begin{theorem}
	\label{TeoSepBq} Let
	$\mathcal{F}_{\omega}$ be the singular holomorphic foliation defined
	by $\omega= H_1 \omega_{ig} + dH_2$ for some $0 \leq i < g$, where $H_j \in \mathbb{C}\{x\}[y]$, $\deg_y H_j<\deg_y F=v_0$ with $j=1,2$, $H_1\neq 0$ and $H_2\in\langle x,y\rangle$. Then the foliation
	$\mathcal{F}_{\omega}$ is dicritical in the last triple point $T_g$ of the
	dual graph $G(\mathcal{C}_F)$ if and only if
	$$I(F,H_1)+v_i+v_g< I(F,H_2).$$  Moreover, all plane branches parameterized by
	$$\psi_{u}(t)=\left(t^{\beta_0},\displaystyle\sum_{\beta_1 \leq j <
			\beta_{g}}c_jt^{j}+ ut^{\beta_{g}}+\displaystyle\sum_{j > \beta_{g}}s_j(u)t^{j}\right),
	$$ with $u\in \mathbb{C}^*$ and $s_j(u)\in \mathbb{C}(u)$ are separatrices of $\mathcal{F}_{\omega}$ whose strict transform intersects transversally the component of the exceptional divisor corresponding to the triple point $T_g$ of $G(\mathcal{C}_F)$.
\end{theorem}
\Dem
Assume that the foliation $\mathcal{F}_\omega$  is dicritical in the last triple point $T_g$ of the dual graph $G(\mathcal{C}_F)$. The separatrices of $\mathcal{F}_\omega$ corresponding to the triple point $T_g$ have parameterizations as given in  (14)  and consequently, they satisfy the same properties as the curve $\mathcal{C}_{F_a}$ described above.
Hence, if $\psi(t)$ is a parametrization of a branch of this dicritical component, then we have that
$$\psi^* \omega\equiv 0.$$
Note that, the first terms which appear in $\psi^* \omega$ are given by
$$\psi^* \omega= k_1t^{I(F,H_1)} (\upsilon_i t^{\upsilon_i} \upsilon_g t^{\upsilon_g -1} -  \upsilon_g t^{\upsilon_g} \upsilon_i t^{\upsilon_i -1}+ \cdots)dt+ (k_2I(F,H_2)t^{I(F,H_2)-1} + \cdots)dt
$$
with $k_1,k_2$ non-zero constants. Hence, if  $I(F,H_2) \leq I(F,H_1)+\upsilon_i+\upsilon_g$, then $\psi^* \omega\not \equiv 0$ against the hypothesis.

Now, assume that $I(F,H_1)+v_i+v_g< I(F,H_2)$. Let $\varphi_a(t)$ be the family of parameterizations  given in
(\ref{family}). Thus $\varphi^{*}_{a}(\omega)=\varphi^*_a(H_1)\varphi^*_a(\omega_{ig})+\varphi^*_a(dH_2)$.
Denoting $u:=a_{\beta_g}$ we will show that it is possible to take  $a_{i}=s_i(u)\in\mathbb{C}(u)$ for every $i>\beta_g$ such that $\varphi_{a}^*(\omega)=0$.

As $deg_yH_j<deg_yF$, by Proposition \ref{Prop-DecompH}, if $H_j\neq 0$ then we can write  $H_j=\sum_{\delta_j}b_{\delta_j}F_0^{\delta_{j0}}\cdot\ldots\cdot F_g^{\delta_{jg}}$ with $b_{\delta_j}\in\mathbb{C}$ and there exist non-negative integers $\gamma_{j0}, \ldots, \gamma_{jg}$ for $j=1,2$ such that $$I(F,H_j)=I(F,F_0^{\gamma_{j0}}\cdot\ldots\cdot F_g^{\gamma_{jg}}).$$

Denoting ${\rm Coeff}(R(t),t^k)$ the coefficient of $t^k$ in $R(t)\in\mathbb{C}\{t\}$, by Lemma \ref{coef-control-Hwij}, for any $k\geq I(F,H_1)+v_i+v_g$  we obtain that
\begin{equation}\label{aqui}
{\rm Coeff}(\varphi^*_a(H_1\omega_{ig}),t^k)=p_{k}(a_{\beta_g},\ldots
,a_{k_{ig}-1})+r_k\cdot a_{\beta_g}^{\gamma_{1g}}\cdot a_{k_{ig}}
\end{equation} with $r_k\in\mathbb{C}^*$ and $k_{ig}:=k-I(F,H_1)-v_i-v_g+\beta_g+1$.

For each $\prod_{l=0}^{g} F_l^{\delta_{2l}}$ in $H_2$, let us denote $$m_{\delta_2}=\max_{0\leq l\leq g}\{l;\ \delta_{2l}\neq 0\}\ \ \ \mbox{and}\ \ \ I_{\delta_2}=I\left (F,\prod_{l=0}^{g} F_l^{\delta_{2l}}\right ).$$

As ${\rm Coeff}(\varphi^*_a(dH_2),t^k)=\sum_{\delta_2}{\rm Coeff}(d(b_{\delta_2}\prod_{l=0}^{g} F_l^{\delta_{2l}}(\varphi_{a})),t^k)$, by Corollary \ref{coef-dH}, for $k\geq I_{\delta_2}$ we get
$${\rm Coeff}(d(b_{\delta_2}\prod_{l=0}^{g} F_l^{\delta_{2l}}(\varphi_{a})),t^k)\in\mathbb{C}[a_{\beta_g},\ldots ,a_{\theta_{\delta_2}}]\ \ \mbox{with}\ \  \theta_{\delta_2}:=k-I_{\delta_2}+\beta_{m_{\delta_2}}+1.$$

But $\beta_{m_{\delta_2}}\leq\beta_g$ and $I_{\delta_2}\geq I(F,H_2)> I(F,H_1)+v_i+v_g$ for any $\delta_2$, so
$$\theta_{\delta_2}=k-I_{\delta_2}+\beta_{m_{\delta_2}}+1 < k-I(F,H_1)-v_i-v_g+\beta_g+1= k_{ig}$$ and, by (\ref{aqui}), we obtain that
$${\rm Coeff}(\varphi^*_a(H_1)\varphi^*_a(\omega_{ig})+\varphi^*_a(dH_2),t^k)=P_k(a_{\beta_g},\ldots ,a_{k_{ig}-1})+r_k \cdot a_{\beta_g}^{\gamma_{1g}} \cdot a_{k_{ig}},$$
for some polynomial $P_k(a_{\beta_g},\ldots ,a_{k_{ig}-1})$ (admitting $H_2=0$).

Remark that $ord_t(\varphi^*_a(\omega))\geq\min\{ord_t(\varphi^*_a(H_1\omega_{ig})), ord_t(\varphi^*_a(dH_2))\}\geq I(F,H_1)+v_i+v_g$. In this way, $\varphi^*_a(\omega)=\varphi^*_a(H_1)\varphi^*_a(\omega_{ig})+\varphi^*_a(dH_2)=0$ is equivalent to solve the system
$$	P_k(a_{\beta_g},\ldots ,a_{k_{ig}-1})+r_k \cdot a_{\beta_g}^{\gamma_{1g}} \cdot a_{k_{ig}}=0
$$ for all $k\geq I(F,H_1)+v_i+v_g$. A such solution exists and it can be obtained by the recurrence relation
\begin{equation}
	\label{Coef-tk-Thm}
	a_{k_{ig}}=-\frac{P_k(a_{\beta_g},\ldots ,a_{k_{ig}-1})}{r_k \cdot a_{\beta_g}^{\gamma_{1g}}},
\end{equation}
since $r_k, a_{\beta_g}\in\mathbb{C}^*$.

In particular, taking $k=I(F,H_1)+v_i+v_g=:k_0$ in the above expression we get $$a_{\beta_g+1}=-\left ( r_{k_0} \cdot a_{\beta_g}^{\gamma_{1g}}\right )^{-1}\cdot P_{k_0}(a_{\beta_g})\in\mathbb{C}(a_{\beta_g})$$ we vanish the coefficient of $t^{k_0}$ in $\varphi^*_a(\omega)$.

Using the previous recurrence relation, we can vanish all terms in $\varphi^*_a(\omega)$ setting the parameters $a_i$ in $\varphi_a(t)$ as a rational function in $\mathbb{C}(a_{\beta_g})$. Hence, considering the parameter $u:=a_{\beta_g}\in\mathbb{C}\setminus\{0\}$ we get the family of parameterizations
$$\psi_{u}(t):=\left(t^{\beta_0},\displaystyle\sum_{\beta_1 \leq j <
	\beta_{g}}c_jt^{j}+ ut^{\beta_{g}}+\displaystyle\sum_{j>\beta_{g}}s_j(u)t^{j}\right),
$$ with $s_j(u):=a_j\in \mathbb{C}(u)$ obtained in (\ref{Coef-tk-Thm}) and satisfying $\psi^*_u(\omega)=0$. As $\psi_u(t)$ defines a family of plane branches with the same characteristic exponents of $\mathcal{C}_F$, every element in the family is topologically equivalent to $\mathcal{C}_F$. This allow us to conclude that the foliation defined by $\omega=0$ is dicritical in the last triple point $T_g$ of the dual graph $G(\mathcal{C}_F)$.
\cqd

If we change $F$ by a semiroot $F_{j+1}$ for $0\leq j< g$ in the previous theorem, then we can describe $1$-forms that define dicritical foliations
in any triple point $T_j$ of the dual graph $G(\mathcal{C}_F)$.

As before, we consider $\{F_0, F_1,\ldots , F_g, F_{g+1}=F\}$ the canonical system of semiroots of $F$, $\Gamma_F=\langle v_0, v_1,\ldots , v_g\rangle$ and $\{\beta_0, \beta_1,\ldots ,\beta_g\}$ the value semigroup and the characteristic exponents of $\mathcal{C}_F$, respectively.

\begin{corollary}
		\label{CorSep} Let
	$\mathcal{F}_{\omega}$ be the singular holomorphic foliation defined
	by $\omega= H_1 \omega_{ij} + dH_2$
	for some $0 \leq i <j\leq g$, where $H_l \in \mathbb{C}\{x\}[y]$, $\deg_y H_l<\deg_y F_{j+1}=\frac{v_0}{e_j}$ with $l=1,2$, $H_1\neq 0$ and $H_2\in\langle x,y\rangle$. The foliation $\mathcal{F}_{\omega}$ is dicritical in the triple point $T_j$ of the dual graph $G(\mathcal{C}_F)$ if and only if
	$$I(F,H_1)+v_i+v_j< I(F,H_2).$$
	Moreover, $\mathcal{F}_{\omega}$ admits a family of separatrices parameterized by
	\begin{equation}
		\label{Sep-family}
		\psi_{u}(t)=\left(t^{\frac{\beta_0}{e_j}},\displaystyle\sum_{\beta_1 \leq l <
			\beta_{j}}c_lt^{\frac{l}{e_j}}+ ut^{\frac{\beta_{j}}{e_j}}+\displaystyle\sum_{l > \frac{\beta_{j}}{e_j}}s_l(u)t^{l}\right),
	\end{equation} with $u\in \mathbb{C}^*$ and $s_l(u)\in \mathbb{C}(u)$.
\end{corollary}
\Dem Notice that $\{F_0,F_1,\ldots ,F_j,F_{j+1}\}$ is the canonical system of semiroots for $F_{j+1}$. So, by (\ref{semiroot-data}), the value semigroup and the characteristic exponents of $\mathcal{C}_{F_{j+1}}$ are respectively,  $\Gamma_{F_{j+1}}=\langle \frac{v_0}{e_j},\ldots ,\frac{v_j}{e_j}\rangle$ and $\left \{\frac{\beta_0}{e_j},\ldots ,\frac{\beta_j}{e_j} \right \}$, where $e_j=GCD(v_0,\ldots ,v_j)$. In addition,
$I(F_{j+1},F_l)=\frac{v_l}{e_{j}}=\frac{I(F,F_l)}{e_j}$ for every $0\leq l\leq j$.

By Proposition \ref{Prop-DecompH} any  $H_l\in\mathbb{C}\{x\}[y]$ with $\deg_yH_l<\deg_yF_{j+1}$ and $l=1,2$ can be expressed by $H_l=\sum_{\delta}b_{\delta_l}F_0^{\delta_{l0}}\cdot\ldots\cdot F_j^{\delta_{lj}}\in\mathbb{C}\{x,y\}$ with $I(F_{j+1},H_l)=I(F_{j+1},F_0^{\gamma_{l0}}\cdot\ldots\cdot F_j^{\gamma_{lj}})$ for some non-negative integers $\gamma_{l0},\ldots ,\gamma_{lj}$. So,
$$
I(F_{j+1},H_l)=\gamma_{l0}\cdot \frac{v_0}{e_j}+\ldots +\gamma_{lj}\cdot \frac{v_j}{e_j}=\frac{\gamma_{l0}\cdot v_0+\ldots +\gamma_{lj}\cdot v_j}{e_j}=\frac{I(F,H_l)}{e_j}.
$$	
Consequently,  $I(F_{j+1},H_1)+\frac{v_i}{e_j}+\frac{v_{j}}{e_j}<I(F_{j+1},H_2)$ if and only if $I(F,H_1)+v_i+v_{j}<I(F,H_2)$.	

Hence, with similar arguments as in the previous theorem, we get that if $\mathcal{F}_{\omega}$ is dicritical in the triple point $T_j$, then the condition $I(F_{j+1},H_1)+\frac{v_i}{e_j}+\frac{v_{j}}{e_j}<I(F_{j+1},H_2)$ must be fulfilled.

Considering the family given by parameterizations
\begin{equation}
	\label{family2}
\varphi_a(t)=\left (t^{\frac{\beta_0}{e_j}},\sum_{\beta_1\leq l<\beta_j}c_lt^\frac{l}{e_j}+\sum_{l \geq \frac{\beta_j}{e_j}}a_lt^l\right )	
\end{equation}
with $a_{\frac{\beta_j}{e_j}}\neq 0$ and proceeding with the same analysis on the coefficients of $\varphi^*_a(\omega)$ as in the previous theorem, in order to obtain $\varphi_a^*(\omega)=0$, for all $k\geq I(F_{j+1},H_1)+\frac{v_i}{e_j}+\frac{v_j}{e_j}$ we must have
\begin{equation}\label{system2}
	P_k\left (a_{\frac{\beta_j}{e_j}},\ldots ,a_{k_{ij}-1}\right )+r_k \cdot a_{\frac{\beta_j}{e_j}}^{\gamma_{1j}} \cdot a_{k_{ij}}=0
\end{equation} where $k_{ij}:=k-I(F_{j+1},H_1)-\frac{v_i}{e_j}-\frac{v_j}{e_j}+\frac{\beta_j}{e_j}+1$ and $r_k\in\mathbb{C}^*$. Hence, we obtain the recurrence relation
$$	a_{k_{ij}}=- \left (r_k \cdot a_{\frac{\beta_j}{e_j}}^{\gamma_{1j}}\right )^{-1}\cdot	P_k\left (a_{\frac{\beta_j}{e_j}},\ldots ,a_{k_{ij}-1}\right ).
$$

In this way, the corollary follows from the previous theorem considering the curve $\mathcal{C}_{F_{j+1}}$, $u:=a_{\frac{\beta_j}{e_j}}\in\mathbb{C}^*$ and $a_l=s_l(u)\in\mathbb{C}(u)$ for $l>\frac{\beta_j}{e_j}$.
\cqd

	Notice that Theorem \ref{TeoSepBq} and Corollary \ref{CorSep} give us a constructive and effective method to present dicritical foliations in a given triple point in the dual graph of a plane branch and to describe parameterizations for the separatrices in such dicritical component up to the desired order.
	
\begin{remark}\label{more}
	Given $\omega=H_1\omega_{ij}+dH_2$ satisfying the hypothesis of the previous corollary and $I(F,H_1)=\sum_{l=0}^{j-1}\gamma_{1l}v_l$, that is, $\gamma_{1j}=0$ then, in (\ref{system2}), we obtain $	P_k(a_{\frac{\beta_j}{e_j}},\ldots ,a_{k_{ij}-1})+r_k \cdot a_{k_{ij}}=0$ and consequently, $a_i$ is a polynomial in $a_{\frac{\beta_j}{e_j}}$ for any $i>\frac{\beta_j}{e_j}$, since $r_k$ is a non-zero constant. In this case, we obtain an extra separatrix for $\mathcal{F}_{\omega}$ taking $a_{\frac{\beta_j}{e_j}}$ with a (not necessarily primitive) parameterization $$\psi_{0}(t)=\left(t^{\frac{\beta_0}{e_j}},\displaystyle\sum_{\beta_1 \leq l <
		\beta_{j}}c_lt^{\frac{l}{e_j}} +\displaystyle\sum_{l>\frac{\beta_{j}}{e_j}}s_l(0)t^{l}\right)
	$$
	and not topologically equivalent to $\mathcal{C}_{F_{j+1}}$.
	
	It is immediate that any irreducible factor $H\in\mathbb{C}\{x,y\}$ of $H_1$ and $H_2$ define a separatrix for $\mathcal{F}_{\omega}$. In addition, if $F_i$ (respectively $F_j$) divides $H_2$, then $F_i$ (respectively $F_j$)
	is a separatrix for $\mathcal{F}_{\omega}$.
\end{remark}

The following examples illustrate the above results.

\begin{example}\label{exemplo-sep} Let us consider the plane branch $\mathcal{C}_F$ with semigroup $\Gamma=\langle 6,9,22\rangle$ as in Example \ref{ex3}. Recall that the characteristics exponents of $\mathcal{C}_F$ are $\beta_0=6, \beta_1=9$ and $\beta_2=13$.
\begin{itemize}
\item Notice that $\zeta_1=(6xy)dy-(9y^2+5x^4)dx=y\cdot\omega_{01}-d(x^5)$ satisfies
$$24=9+6+9=I(F,H_1)+v_0+v_1<I(F,H_2)=30,$$
consequently by Corollary \ref{CorSep}, $\mathcal{F}_{\zeta_1}$ is dicritical admitting a family of separatrices in the first triple point of $G(\mathcal{C}_F)$ parameterized by
$$\psi_u(t)=\left(t^2, \ ut^3+\frac{5}{6u}t^{5}-\frac{25}{72u^3}t^{7}+\frac{125}{432u^5}t^9+\sum_{i\geq 11} q_1(u)t^i\right)$$
with $q_1(u)\in\mathbb{C}[u^{-1}]$. Moreover, by Remark \ref{monomial}, $\mathcal{F}_{\zeta_1}$ also admits the separatrices $(0,t)$, $\left (t, \frac{\sqrt{15}}{3}t^2\right )$ and $\left (t, -\frac{\sqrt{15}}{3}t^2\right )$.

\item Taking $\zeta_2=x(12y^2-12x^2y+x^5)dy+2y(2x^3+10x^2y+x^4-11y^2+3x^5)dx=y\cdot \omega_{02}+d(x^6y)$ we have
$$37=9+6+22=I(F,H_1)+v_0+v_2<I(F,H_2)=45.$$
So, by Theorem \ref{TeoSepBq}, $\mathcal{F}_{\zeta_2}$ is a dicritical foliation in the last triple point of $G(\mathcal{C}_F)$. Moreover, the family
$$\psi_u(t)=\left (t^6, t^9+t^{12}+ut^{13}-\frac{u^2}{2}t^{17}+\left (-\frac{15}{32}+\frac{u^3}{2}\right )t^{21}-\frac{1}{44}t^{24}+\sum_{i\geq 25}q_2(u)t^i\right )$$
with $q_2(u)\in\mathbb{C}[u]$ describe separatrices for $\mathcal{F}_{\zeta_2}$. Others separatrices for $\mathcal{F}_{\zeta_2}$ are $(0,t)$, $(t,0)$
and, by Remark \ref{more}, $\psi_0(t)=\left (t^6,t^9+t^{12}-\frac{15}{32}t^{21}-\frac{1}{44}t^{24}+\sum_{i\geq 9} q_2(0)t^{3i}\right )$, that is,
$$\left (t^2,t^3+t^{4}-\frac{15}{32}t^7-\frac{1}{44}t^8+\sum_{i\geq 9} q_2(0)t^{i}\right ).$$

\item Considering $$\hspace{-1.5cm}\zeta_3=\left (2x(11x^3-11x^4-2y^2)+y\left (\frac{227}{10}x^3+\frac{33}{5}y^2-\frac{99}{10}x^2y+\frac{33}{10}x^4\right )\right )dy+$$
$$\hspace{4cm} +xy\left (\frac{33}{20}y+9x\right )(-3x-4y+4x^2 )dx$$ we can write $\zeta_3=x\cdot\omega_{12}+d\left (\frac{33}{20}y^2F_2\right )$
with $F_2=y^2-x^3-2x^2y+x^4$. As
$$37=6+9+22=I(F,H_1)+v_1+v_2<I(F,H_2)=2\cdot 9+22=40$$
the previous results ensure that $\mathcal{F}_{\zeta_3}$ is dicritical in the last triple point of $G(\mathcal{C}_F)$ and
$$\psi_u(t)=\left ( t^6,t^9+t^{12}-ut^{13}+\frac{35}{18}u^2t^{17}+\frac{473}{180}ut^{19}-\frac{748}{189}u^2t^{20}+\sum_{i\geq 21}q_3(u)t^{i}\right )$$ define separatrices for $\mathcal{F}_{\eta_3}$. By Remark \ref{more}, the curves $(t,0)$ and $\psi_0(t)=(t^2,t^3+t^4)$ (that is, the curve defined by  $F_2=y^2-x^3-2x^2y+x^4$) are also separatrices for $\mathcal{F}_{\zeta_3}$.
\end{itemize}	
\end{example}
Let $\mathcal{C}_F$ be the plane branch with semigroup $\Gamma=\langle 6,9,22\rangle$ as in Example \ref{ex3} and $$\zeta=(6xy)dy-(9y^2+4x^3)dx=y\cdot\omega_{01}-d(x^4).$$ In this case $24=9+6+9=I(F,H_1)+v_0+v_1=I(F,H_2)$ and the foliation $\mathcal{F}_{\zeta}$ is not dicritical.  The unique separatrix of the foliation $\mathcal{F}_{\zeta}$ is the curve $x=0$. Note that $\mathcal{F}_{\zeta}$ is not a second type foliation: there is a saddle-node singularity in one of the corners of its reduction of singularities (see \cite{Mat-S}).

\section{Analytical invariants of $\mathcal{C}_F$ and dicritical foliations}
\label{Section-Invariants}

As before, $\mathcal{C}_F$ is a plane branch defined by a Weierstrass polynomial $F\in\mathbb{C}\{x\}[y]$  with $mult(F)=v_0$ admitting a parameterization $\varphi(t)=(x(t),y(t))$. Considering $\varphi^*:\mathbb{C}\{x,y\}\rightarrow\mathbb{C}\{t\}$ defined by
$\varphi^*(H)=H(x(t),y(t))$ we have that $\ker\varphi^*=\langle F\rangle$
and $Im\varphi^*=\mathbb{C}\{x(t),y(t)\}$. So, we obtain the exact
sequence of $\mathbb{C}$-algebras
$$\{0\}\rightarrow\langle F\rangle\hookrightarrow\mathbb{C}\{x,y\}\rightarrow\mathbb{C}\{x(t),y(t)\}\rightarrow\{0\}.$$
In this way, the local ring of $\mathcal{C}_F$ is
$\mathcal{O}:=\frac{\mathbb{C}\{x,y\}}{\langle
	F\rangle}\cong\mathbb{C}\{x(t),y(t)\}\subset
\mathbb{C}\{t\}=:\overline{\mathcal{O}}$ where
$\overline{\mathcal{O}}$ denotes the integral closure of
$\mathcal{O}$ in its field of fractions.

Let $\nu:\overline{\mathcal{O}}\rightarrow\mathbb{Z}_{\geq
	0}\cup\{\infty\}$ be the discrete normalized valuation given by
$\nu(p(t))=ord_t p(t)$ for $p(t)\in\mathbb{C}\{t\}$ ($\nu(0)=\infty$) and denote
$\nu(H)=\nu(\varphi^*(H))$ for $H\in\mathbb{C}\{x,y\}$. In this way,
the value semigroup of $\mathcal{C}_F$ is given by $\Gamma_F=\nu(\mathcal{O})$. In addition, the conductor ideal  $(\mathcal{O}:\overline{\mathcal{O}}):=\{h\in\mathcal{O};\ h\overline{\mathcal{O}}\subseteq\mathcal{O}\}$ of $\mathcal{O}$ in $\overline{\mathcal{O}}$ satisfies $(\mathcal{O}:\overline{\mathcal{O}})=\langle t^{\mu_F}\rangle$, that is, if $p(t)\in\mathbb{C}\{t\}$ is such that $\nu(p(t))\geq\mu_F$, then there exists $H\in\mathbb{C}\{x,y\}$ with $\varphi^*(H)=p(t)$. The integer $\mu_F$ is called the conductor of $\Gamma_{F}$.

If $H\in\langle x,y\rangle$ by (\ref{phi-omega}) we get $ord_t(\varphi^*(dH))=\nu(H)-1$, that is, $ord_t(\varphi^*(dH))+1\in\Gamma_{F}$.
In this way, given $\omega\in\Omega^1$ such that $\varphi^*(\omega)\neq 0$
we define the value of $\omega$ as $\nu(\omega):=ord_t(\varphi^*(\omega))+1$. Setting $\nu(\omega)=\infty$ if $\varphi^*(\omega)=0$, we define
$$\Lambda_F:=\left \{\nu(\omega);\ \omega\in\Omega^1\right
\}\supseteq\Gamma_F\setminus\{0\}.$$

\begin{remark}\label{elimination}
The set $\Lambda_F\subset\overline{\mathbb{N}}$ is an analytic invariant for $\mathcal{C}_F$ and it is the main ingredient for the analytic classification of plane branches presented in \cite{HH} (see \cite{handbook} for an extended version).

In particular, the set $\Lambda_{F}$ allow us to identify terms in a parameterization of $\mathcal{C}_F$ that can be eliminated by change of parameter and coordinates. More specifically, by Proposition 1.3.11 and Theorem 1.3.9 in \cite{handbook}, given a plane branch $\mathcal{C}_F$ with Puiseux parameterization $\left (t^{v_0},\sum_{i> v_0}a_it^i \right )$ if there exists $\omega=Adx+Bdy\in\Omega^1$ with $\nu(\omega)=k+v_0$, $A\in\langle x,y\rangle^2$ and $B\in\langle x^2,y\rangle$ then $\mathcal{C}_F$ is analytically equivalent to a plane branch with parameterization
$\left (t^{v_0},\sum_{i> v_0}b_it^i \right )$ where $b_i=a_i$ for $i<k$ and $b_k=0$.
\end{remark}

We can define $\Lambda_F$ by means the $\mathcal{O}$-module of K\"{a}hler differentials of ${\mathcal	O}$ (or $\mathcal{C}_F$), that is,
$$\Omega_{\mathcal{O}}:=\Omega_{{\mathcal		O}/\mathbb{C}}=\frac{{\mathcal	O}dx+{\mathcal O}dy}{{\mathcal O}(F_xdx+F_ydy)}\cong \frac{\Omega^{1}}{\mathcal{F}(\mathcal{C}_F)},$$ where $\mathcal{F}(\mathcal{C}_F):=F\cdot\Omega^{1}+\mathbb{C}\{x,y\}\cdot dF$.

If $\eta\in\mathcal{F}(\mathcal{C}_F)$ then $\varphi^*(\eta)=0$ and $\varphi^*(\omega+\eta)=\varphi^*(\omega)$ for any $\omega\in\Omega^1$. Thus, given $\overline{\omega}=\omega+\mathcal{F}(\mathcal{C}_F)\in\Omega_{\mathcal{O}}$ we can define $\varphi^*(\overline{\omega}):=\varphi^*(\omega)$ and $\nu(\overline{\omega}):=\nu(\omega)$. For any singular plane branch, the torsion submodule $\mathcal{T}:=\{\overline{\omega}\in\Omega_{\mathcal O};\ h\overline{\omega}=0\ \mbox{for some}\ h\in\mathcal{O}\setminus\{0\}\}\subset\Omega_{\mathcal O}$ is non trivial and we can rewrite $\mathcal{T}=\{\overline{\omega}\in\Omega_{\mathcal O};\ \varphi^*(\overline{\omega})=0\}$. In particular, we have $\frac{\Omega_{\mathcal O}}{\mathcal{T}}\cong\varphi^*(\Omega_{\mathcal O})=\varphi^*(\Omega^1)\subset\mathbb{C}\{t\}$
and $\Lambda_F=\{\nu(\overline{\omega});\ \overline{\omega}\in\Omega_{\mathcal O}\setminus\mathcal{T}\}$ (see Section 7.1 in \cite{coloquio} and \cite{basestandard}).

There exists a finite subset
$L=\{\ell_1,\ldots ,\ell_k\}\subset \Lambda_F$ such that any
$\ell\in\Lambda_F$ can be expressed as
$\ell=\ell_i+\gamma$ for some $\gamma\in\Gamma_F$ and
$\ell_i\in L$, that is, the set $\Lambda_F$ is a finitely generated
$\Gamma_F$-monomodule. A set $\mathcal{G}=\{\overline{\omega_1},\ldots ,\overline{\omega_k}\}\subset \frac{\Omega_{\mathcal O}}{\mathcal{T}}$ such that
$\nu(\overline{\omega_i})=\nu(\omega_i)=\ell_i\in L$ is  a set of generators for $\frac{\Omega_{\mathcal O}}{\mathcal{T}}$ as $\mathcal{O}$-module and it is called a {\bf Standard Basis} of $\frac{\Omega_{\mathcal O}}{\mathcal{T}}$.

Fixing a minimal set of generators $L$ for
$\Lambda_F$, that is, $L$ is a set of generators for $\Lambda_F$ (as
$\Gamma_F$-monomodule) and $\ell_i\not\in\ell_j+\Gamma_F$ for
$\ell_i,\ell_j\in L$ and $i\neq j$, we call $\mathcal{G}=\{\overline{\omega_i};\ \nu(\omega_i)\in L\}$ a {\bf Minimal Standard Basis} for $\frac{\Omega_{\mathcal O}}{\mathcal{T}}$.
In \cite{basestandard} we provide an algorithm to compute a (minimal) Standard Basis $\mathcal{G}$ for $\frac{\Omega_{{\mathcal O}}}{\mathcal{T}}$ by means a parameterization $\varphi(t)$ and in Section 7.3 of \cite{coloquio} we describe a method to obtain $\mathcal{G}$ using $F$.

\begin{remark}\label{jac}
The set $\Lambda_F$ determines and it is determined by the values of elements in the Jacobian ideal
$J_F:=\langle F_x,F_y\rangle$ in $\mathcal{O}$. In
fact, we have the isomorphism (as $\mathcal{O}$-module)
$$\begin{array}{cccc} \Psi : & \mathcal{O} F_y+\mathcal{O} F_x &	\rightarrow &
	\varphi^*(\Omega_{{\mathcal O}})\cong \frac{\Omega_{\mathcal{O}}}{\mathcal{T}} \\
	& AF_y+BF_x & \mapsto & \varphi^*(Adx-Bdy).
\end{array}$$

Notice that $0=dF=F_xdx+F_ydy\in\Omega_{{\mathcal O}}$ then given $\omega=Adx-Bdy\in\Omega_{{\mathcal O}}$ we have that
$$F_y\omega=(AF_y+BF_x)dx-BdF=(AF_y+BF_x)dx \ \ (\mbox{in}\ \Omega_{{\mathcal O}}).$$ As $F\in\mathbb{C}\{x\}[y]$ is a Weierstrass polynomial with $mult(F)=v_0$, then $\nu(F_y)=\mu_F+v_0-1$ (see Corollary 7.16 in \cite{hefez}) and $\nu(dx)=v_0$. Thus, by the above expression, we get
$\nu(AF_y+BF_x)=\nu(F_y \omega)-v_0=\mu_F-1+\nu(\omega).$
In this way, $$\nu(J_F)=\nu(\mathcal{O} F_x+\mathcal{O} F_y)=\mu_F-1+\Lambda_F.$$
Pol, in \cite{pol}, generalizes this result for a reduced complete intersection  curve.
\end{remark}

Notice that if $\ell$ is an element of a minimal set of generators $L$ for $\Lambda_F$, then any $\omega=Adx+Bdy\in\Omega^1$ such that $\nu(\omega)=\ell$ defines a foliation, that is, $GCD(A,B)=1$. Indeed, if $A=A_1G$ and $B=B_1G$ with $G\in\langle x,y\rangle$ then $\ell=\nu(\omega)=\nu(G)+\nu(A_1dx+B_1dy)$ that is a contradiction by minimality of $L$.

Let $\omega=P_1dx+P_2dy\in\Omega^1$. By the Weierstrass Division Theorem we can express $P_2=Q_2F_y+A_2$ and $P_1-Q_2F_x=Q_1F+A_1$, where $Q_1, Q_2\in\mathbb{C}\{x,y\}$, $A_1, A_2\in\mathbb{C}\{x\}[y]$, $deg_y(A_1)<deg_y(F)=v_0$ and $deg_y(A_2)<deg_y(F_y)=v_0-1$. Thus
$$\omega=P_1dx+P_2dy=A_1dx+A_2dy+FQ_1dx+Q_2dF.$$
As $FQ_1dx+Q_2dF\in\mathcal{F}(\mathcal{C}_F)$ we get $\nu(\omega)=\nu(A_1dx+A_2dy)$ and consequently
\begin{equation}\label{economy}
\Lambda_F=\{\nu(Adx+Bdy);\ A, B\in\mathbb{C}\{x\}[y]\ \mbox{with}\ deg_y(A)<v_0\ \mbox{and}\ deg_y(B)<v_0-1\}.
\end{equation}

Using Corollary \ref{CorSep} we can relate elements of $\Lambda_F$ with the contact of $\mathcal{C}_F$ and a curve corresponding to a particular element of the family $\psi_{u}(t)$ as in (\ref{Sep-family}).

For commodity to the reader we recall some results concerning the contact and the intersection multiplicity of plane curves (see \cite{hefez} or \cite{wall}).

Consider two irreducible plane curves $\mathcal{C}_F$ and $\mathcal{C}_G$ with Puiseux parameterizations given respectively by $(t^{v_0},\phi(t))$ and $(t^{v'_0},\phi'(t))$. The contact $c(\mathcal{C}_F,\mathcal{C}_G)$ of $\mathcal{C}_F$ and $\mathcal{C}_G$ is defined by
$$c(\mathcal{C}_F,\mathcal{C}_G)=\max_{\gamma,\delta}\frac{ord_t(\phi(\gamma t^{v'_0})-\phi'(\delta t^{v_0}))}{v_0v'_0}$$
where $\gamma,\delta\in\mathbb{C}$ with $\gamma^{v_0}=1=\delta^{v'_0}$.

In what follows we take parameterizations of $\mathcal{C}_F$ and $\mathcal{C}_G$ such that the maximum in the previous expression is achieved.

Remark that, by definition, the series $\phi(t^{v'_0})$ and $\phi'(t^{v_0})$ coincide up to the order $c(\mathcal{C}_F,\mathcal{C}_G)v_0v'_0-1$. In addition, for any irreducible plane curve $\mathcal{C}_H$ we have that
$$c(\mathcal{C}_F,\mathcal{C}_H)\geq\min\{c(\mathcal{C}_F,\mathcal{C}_G),c(\mathcal{C}_H,\mathcal{C}_G)\}$$
and the two smallest numbers among these three coincide.

Let $\beta_i, e_i, n_i$ and $v_i$ be the integers defined at the beginning of Section 2, related to the branch $\mathcal{C}_F$. We indicate by $\beta'_i, e'_i, n'_i$ and $v'_i$ the respective integers for $\mathcal{C}_G$.

The contact of two branches $\mathcal{C}_F$ and $\mathcal{C}_G$ is related to the intersection multiplicity $I(F,G)$ in the following way (see \cite{Merle}):

If $c(\mathcal{C}_F,\mathcal{C}_G)<\frac{\beta_1}{v_0}$, then $I(F,G)=c(\mathcal{C}_F,\mathcal{C}_G)v_0v_0'$. Moreover  $\frac{\beta_q}{v_0}\leq c(\mathcal{C}_F,\mathcal{C}_G) <\frac{\beta_{q+1}}{v_0}$ for some $q \in \{1, \ldots, g\}$ if and only if
\begin{equation}\label{inter}\frac{I(F,G)}{v'_0}=\frac{n_qv_q+v_0 \cdot c(\mathcal{C}_F,\mathcal{C}_G)-\beta_q}{n_0\cdot\ldots\cdot n_q}.\end{equation}

Let us consider $\omega=H_1 \omega_{ij}+dH_2$ with $0\leq i<j\leq g$ satisfying the hypothesis of Corollary \ref{CorSep}, that is,
$\omega$ defines a dicritical foliation in the $j$th triple point of the dual graph of $\mathcal{C}_F$ and $\mathcal{F}_{\omega}$ admits a family of separatrices parameterized by $\psi_u(t)$ as (\ref{Sep-family}).

Let $F_u\in\mathbb{C}\{x\}[y]$ be the irreducible Weierstrass polynomial such that $F_u(\psi_u(t))=0$. In particular, $\mathcal{C}_{F_u}$ and $\mathcal{C}_{F_{j+1}}$ are topologically equivalent and consequently they admit the same characteristic exponents $\left \{ \beta'_0:=\frac{\beta_0}{e_j},\ldots ,\beta'_j:=\frac{\beta_j}{e_j}\right \}$ and the same value semigroup $\langle v'_0:=\frac{v_0}{e_j},\ldots ,v'_j:=\frac{v_j}{e_j}\rangle$.

\begin{remark}\label{contact}
	Notice that  $\frac{\beta'_{j}}{v'_0}=\frac{\beta_{j}}{v_0}\leq c(\mathcal{C}_{F},\mathcal{C}_{F_u})\leq \frac{\beta_{j+1}}{v_0}$ for any $u\in\mathbb{C}^*$ and, by (\ref{Sep-family}), $\frac{\beta'_{j}}{v'_0}=\frac{\beta_{j}}{v_0}= c(\mathcal{C}_{F},\mathcal{C}_{F_u})$ if and only if $u\neq c_{\beta_{j}}$. In fact, if $\frac{\beta_{j+1}}{v_0}<c(\mathcal{C}_{F},\mathcal{C}_{F_u})$ then in $\psi_u(t)$ we should have a term with exponent $\frac{\beta_{j+1}v'_0}{v_0}=\frac{\beta_{j+1}}{e_j}\not\in\mathbb{Z}_{\geq 0}$, that is an absurd.
\end{remark}

By above remark, for any $ u\neq c_{\beta_{j}}$ we conclude, by (\ref{inter}), that
$$I(F,F_u)=\frac{v'_0n_jv_{j}}{n_0\cdot\ldots\cdot n_{j}}=\frac{v'_0n_{j}v_{j}}{\frac{v_0}{e_{j}}}=\frac{v'_0n_{j}v_{j}}{v'_0}=n_{j}v_{j}.$$

So, for $u=c_{\beta_j}$ the curve $\mathcal{C}_{F_u}$ has a special behavior. In order to simplify the notations we denote
\begin{equation}\label{star}\psi_{\star}(t):=\psi_{c_{\beta_j}}(t)\ \ \ \mbox{and}\ \ \ F_{\star}:=F_{c_{\beta_j}}.\end{equation}

If $F$ is not a separatrix of $\omega$, that is, $\infty\neq\nu(\omega)\in\Lambda_F$, then we can relate $c(\mathcal{C}_F,\mathcal{C}_{F_{\star}})$ and $\nu(\omega)$.

\begin{theorem}\label{cont-value} Given $\omega=H_1 \omega_{ij}+dH_2$ with $0\leq i<j\leq g$ such that $I(F,H_1)+v_i+v_{j}< I(F,H_2)$, where $H_l \in \mathbb{C}\{x\}[y]$, $\deg_y H_l<\deg_y F_{j+1}=\frac{v_0}{e_j}$ for $l=1,2$, $H_1\neq 0$ and $H_2\in\langle x,y\rangle$
	we have
	$$c(\mathcal{C}_{F},\mathcal{C}_{F_{\star}})=\frac{\nu(\omega)-I(F,H_1)-v_i-v_j+\beta_j}{v_0}.$$	
\end{theorem}
\Dem Given $\omega=A(x,y)dx+B(x,y)dy\in\Omega^1$, $\phi(t)=(x(t),y(t))\in\mathbb{C}\{t\}\times\mathbb{C}\{t\}$ and $n\in\mathbb{Z}_{>0}$ we denote
$$(\omega\phi)(t)=A(x(t),y(t))x'(t)+B(x(t),y(t))y'(t)=\phi^*(\omega);$$
$$(\omega\phi)(t^n)=A(x(t^n),y(t^n))x'(t^n)+B(x(t^n),y(t^n))y'(t^n);$$
$$\omega(\phi(t^n))=A(x(t^n),y(t^n))(x(t^n))'+B(x(t^n),y(t^n))(y(t^n))'=nt^{n-1}(\omega\phi)(t^n).$$
Consequently,
\begin{equation}\label{tn}
	\begin{array}{c}
	{\rm Coeff}((\omega\phi)(t^{n}),t^{kn})=	{\rm Coeff}((\omega\phi)(t),t^{k})\\ nt^{n-1}{\rm Coeff}((\omega\phi)(t^n),t^{kn})=	{\rm Coeff}(\omega(\phi(t^n)),t^{n(k+1)-1}).
	\end{array}
\end{equation}

By the proof of Corollary \ref{CorSep}, for any member of the family $\varphi_a(t)$, given in (\ref{family2}), we obtain that ${\rm Coeff}((\omega\varphi_a)(t),t^k)\in\mathbb{C}[a_{\beta'_j},\ldots ,a_{\epsilon}]$ where $\epsilon=k-I(F_{j+1},H_1)+\beta'_j-v'_j-v'_i+1$. So, by (\ref{tn}), the coefficients of terms with order up to $\epsilon$ in $\varphi_a(t)$, or equivalently, the coefficients of terms with order up to $v_0\epsilon$ in $\varphi_a(t^{v_0})$ determine all coefficients of terms with order up to
$$\begin{array}{ll}
\epsilon+I(F_{j+1},H_1)-\beta'_j+v'_j+v'_i-1 & \mbox{in}\ (\omega\varphi_a)(t);\\
v_0(\epsilon+I(F_{j+1},H_1)-\beta'_j+v'_j+v'_i-1) & \mbox{in}\ (\omega\varphi_a)(t^{v_0});\\
v_0(\epsilon+I(F_{j+1},H_1)-\beta'_j+v'_j+v'_i)-1 & \mbox{in}\ \omega(\varphi_a(t^{v_0})).
\end{array}$$

As $\psi_{\star}(t^{v_0})$ is a member of the family $\varphi_a(t^{v_0})$ and
${\rm Coeff}((\omega\psi_{\star})(t),t^k)=0$ for all $k$, by (\ref{tn}), we have that ${\rm Coeff}(\omega(\psi_{\star}(t^{v_0}),t^k)=0$ for all $k$. But  $\psi_{\star}(t^{v_0})$ and $\varphi(t^{v'_0})$ coincide up to the order $v_0v'_0c(\mathcal{C}_F,\mathcal{C}_{F_{\star}})-1$ and $0\neq{\rm Coeff}(\psi_{\star}(t^{v_0}),t^{v_0v'_0c(\mathcal{C}_F,\mathcal{C}_{F_{\star}})})\neq {\rm Coeff}(\varphi(t^{v'_0}),t^{v_0v'_0c(\mathcal{C}_F,\mathcal{C}_{F_{\star}})})$, then
$$0={\rm Coeff}(\omega(\psi_{\star}(t^{v_0})),t^k)={\rm Coeff}(\omega(\varphi(t^{v'_0})),t^k)$$
for all $$k<v_0(v'_0c(\mathcal{C}_F,\mathcal{C}_{F_{\star}})+I(F_{j+1},H_1)-\beta'_{j}+v'_j+v'_i)-1=v_0v'_0c(\mathcal{C}_F,\mathcal{C}_{F_{\star}})+v'_0(I(F,H_1)-\beta_{j}+v_j+v_i)-1,$$
recall that $v'_l=\frac{v_l}{e_j}$,  $\beta'_l=\frac{\beta_l}{e_j}$ for $0\leq l\leq j$ and $\frac{I(F,H_1)}{e_j}=I(F_{j+1},H_1)$.

Moreover, we get ${\rm Coeff}(\omega(\varphi(t^{v'_0})),t^k)\neq 0$
for $k=v_0v'_0c(\mathcal{C}_F,\mathcal{C}_{F_{\star}})+v'_0(I(F,H_1)-\beta_{j}+v_j+v_i)-1$. So, by (\ref{tn}), $k=v'_0(ord_t((\omega\varphi)(t)))+v'_0-1$, that is, $k=v'_0(ord_t((\omega\varphi)(t))+1)-1=v'_0\nu(\omega)-1$.

In this way,
$v'_0\nu(\omega)-1=k=v_0v'_0c(\mathcal{C}_F,\mathcal{C}_{F_{\star}})+v'_0(I(F,H_1)-\beta_{j}+v_j+v_i)-1,$
that is,
$$c(\mathcal{C}_F,\mathcal{C}_{F_{\star}})=\frac{\nu(\omega)-I(F,H_1)-v_i-v_j+\beta_{j}}{v_0}.$$
\cqd

Let us illustrate the previous theorem.

\begin{example}\label{exemplo-value} Let $\mathcal{C}_F$ be a plane branch  with semigroup $\Gamma=\langle 6,9,22\rangle$ and Puiseux parameterization $\varphi(t)=(t^6,t^9+t^{12}+2t^{13})$ as in Example \ref{exemplo-sep}.

Considering $\zeta_1=y\cdot\omega_{01}-d(x^5)$ we have $\nu(\zeta_1)=27$. In this case $u=c_{\beta_1}=1$, $$\psi_{\star}(t)=\left(t^2, \ t^3+\frac{5}{6}t^{5}-\frac{25}{72}t^{7}+ h.o.t.\right )\ \ \ \mbox{and}\ \ \
c(\mathcal{C}_{F},\mathcal{C}_{F_{\star}})=\frac{\nu(\zeta_1)-I(F,y)-v_0-v_1+\beta_1}{6}=2.$$
		
For $\zeta_2=y\cdot \omega_{02}+d(x^6y)$ we have $\nu(\zeta_2)=41$, $u=c_{\beta_2}=2$,
$$\psi_{\star}(t)=\left (t^6, t^9+t^{12}+2t^{13}-2t^{17}+h.o.t.\right )\  \ \mbox{and}\ \ c(\mathcal{C}_{F},\mathcal{C}_{F_{\star}})=\frac{\nu(\zeta_2)-I(F,y)-v_0-v_2+\beta_2}{6}=\frac{17}{6}.$$
			
Given $\zeta_3=x\cdot\omega_{12}+d\left (\frac{33}{20}y^2F_2\right )$ we get $\nu(\zeta_3)=41$, $u=c_{\beta_2}=2$,
		$$\psi_{\star}(t)=\left ( t^6,t^9+t^{12}+2t^{13}+\frac{70}{9}t^{17}+h.o.t.\right )\ \  \mbox{and}\ \ c(\mathcal{C}_{F},\mathcal{C}_{F_{\star}})=\frac{\nu(\zeta_3)-I(F,x)-v_1-v_2+\beta_2}{6}=\frac{17}{6}.$$
\end{example}

By (\ref{inter}) we can determine the intersection multiplicity $I(F,F_{\star})$ by means the contact $c(\mathcal{C}_{F},\mathcal{C}_{F_{\star}})$ and consequently, by Theorem \ref{cont-value}, we can relate $I(F,F_{\star})$ and $\nu(\omega)$.

\begin{corollary}
	\label{Inter}
With the previous notations, we have $$I(F,F_{\star})= \nu(\omega)-I(F,H_1)+(n_j-1)v_j-v_i.$$
In particular,
if $g=1$ then $\nu(\omega)=I(F,H_1\cdot F_{\star})-(\mu_F-1)$.
\end{corollary}
\Dem By Remark \ref{contact}, we have $\frac{\beta_j}{v_0}< c(\mathcal{C}_{F},\mathcal{C}_{F_{\star}})\leq \frac{\beta_{j+1}}{v_0}$.

If $c(\mathcal{C}_{F},\mathcal{C}_{F_{\star}})\neq \frac{\beta_{j+1}}{v_0}$, by (\ref{inter}), we get
$$I(F,F_{\star})=v'_0\cdot\left (\frac{n_jv_j+v_0 \cdot c(\mathcal{C}_{F},\mathcal{C}_{F_{\star}})-\beta_j}{n_0\cdot\ldots\cdot n_j}\right )=\nu(\omega)-I(F,H_1)+(n_j-1)v_j-v_i,$$ where the last equality is obtained  by Theorem \ref{cont-value} remembering that $n_0\cdot\ldots\cdot n_j=\frac{v_0}{e_j}=v'_0$.

If $c(\mathcal{C}_{F},\mathcal{C}_{F_{\star}})= \frac{\beta_{j+1}}{v_0}$, then by previous theorem we have
$\beta_{j+1}=\nu(\omega)-I(F,H_1)-v_i-v_j+\beta_j$ that is, by (\ref{relations0}), $\nu(\omega)-I(F,H_1)+(n_j-1)v_j-v_i=v_{j+1}$. On the other hand, by (\ref{inter}), we obtain
$$I(F,F_{\star})=v'_0\cdot\left (\frac{n_{j+1}v_{j+1}+\beta_{j+1}-\beta_{j+1}}{n_0\cdot\ldots\cdot n_{j+1}}\right )=v_{j+1}$$ that gives us the result.

In particular, if $g=1$, we have $i=0$ and $j=g=1$. So, by (\ref{milnor}) we have $\nu(\omega)=I(F,H_1\cdot F_{\star})-(n_1-1)v_1+v_0=I(F,H_1\cdot F_{\star})-(\mu_F-1)$.
\cqd

If the value semigroup of $\mathcal{C}_F$ is $\Gamma_F=\langle v_0,v_1\rangle$, that is, $g=1$ then Lemma \ref{azevedo} ensures that any $1$-form $\eta=Adx+Bdy\in\Omega^1$ can be written as $\eta=P_1 \omega_{01}+dP_2$ with $P_1, P_2\in\mathbb{C}\{x,y\}$. Considering $P_i=Q_iF+H_i$ such that $Q_i\in\mathbb{C}\{x,y\}$ with $H_i\in\mathbb{C}\{x\}[y]$ and $deg_y(H_i)<deg_y(F)=v_0$ we have $\eta=H_1\omega_{01}+dH_2+F\cdot (Q_1\omega_{01}+dQ_2)+Q_2dF$.

In this way, for any  $\eta\in\Omega^1$ there exists $\omega=H_1\omega_{01}+dH_2$  with $H_i\in\mathbb{C}\{x\}[y]$ and $deg_y(H_i)<deg_y(F)=v_0$ for $i= 1, 2$ such that $\nu(\eta)=\nu(\omega)$.
Moreover, if $\nu(\omega)\in\Lambda_F\setminus\Gamma_F$ then  $$I(F,H_2)=\nu(dH_2)\geq\nu(H_1\cdot\omega_{01})=I(F,H_1)+\nu(\omega_{01})>I(F,H_1)+v_1+v_0.$$
that is, if $\nu(\omega)\in\Lambda_F\setminus\Gamma_F$ and $\omega$ defines a foliation, then $\omega$ satisfies the hypothesis of Theorem \ref{TeoSepBq}.

\begin{corollary}\label{Lambda-g-1} Let $\mathcal{C}_{F}$ be a plane branch with value semigroup $\Gamma_{F}=\langle v_0,v_1\rangle$. If an element $\ell\in\Lambda_F\setminus\Gamma_F$ belongs to a minimal set of generators for $\Lambda_F$, then there exists a dicritical foliation $\mathcal{F}_{\omega}$ in the triple point of $G(\mathcal{C}_F)$ defined by $\omega=H_1\cdot\omega_{01}+dH_2\in\Omega^1$ with $\nu(\omega)=\ell$. In particular,  $$\nu(\omega)+\mu_F-1=I(F,H_1\cdot F_{\star})=I\left (F,\frac{\omega\wedge dF}{dx\wedge dy}\right ),$$
	where $F_{\star}$ is defined in (\ref{star}).
\end{corollary}
\Dem By previous comments, if $\ell\in\Lambda_F\setminus\Gamma_F$ belongs to a minimal set of generators for $\Lambda_F$, then there exists $\omega\in\Omega^1$ with $\nu(\omega)=\ell$, satisfying the hypothesis of Theorem \ref{TeoSepBq}, that is, $\mathcal{F}_{\omega}$ defines a dicritical foliation in the triple point of $G(\mathcal{C}_f)$. In particular, by Corollary \ref{Inter} and Remark \ref{jac}, we get $\nu(\omega)+\mu_F-1=I(F,H_1\cdot F_{\star})=I\left (F,\frac{\omega\wedge dF}{dx\wedge dy}\right )$.
\cqd

The above result guarantees that for plane branches with value semigroup $\Gamma_{F}=\langle v_0,v_1\rangle$, the minimal generators for $\Lambda_F$, distinct of $v_0$ and $v_1$, can be obtained considering dicritical foliations in the triple point of $G(\mathcal{C}_F)$.

The previous corollary was obtained by Cano, Corral and Senovilla-Sanz in \cite{nuria} by other methods.

\subsection{The Zariski invariant of $\mathcal{C}_F$}

Now we return to the general case, that is, plane branches $\mathcal{C}_F$ with value semigroup $\Gamma_F=\langle v_0,\ldots ,v_g\rangle$ with $g\geq 1$. Without loss of generality, we can assume that $\mathcal{C}_F$ is defined by a Weierstrass polynomial $F\in \mathbb{C}\{x\}[y]$ satisfying $e_1=GCD(v_0,v_1)<v_0=I(F,x)<v_1=I(F,y)$.

In \cite{torsion}, Zariski shows that $\lambda:=\min(\Lambda_{F}\setminus\Gamma_{F})-v_0$ can be computed directly by a Puiseux parameterization $(t^{\beta_0},\sum_{i\geq \beta_1}c_it^i)$ of $\mathcal{C}_F$. More precisely,
$$\lambda=\min\{i;\ c_i\neq 0\ \mbox{and}\ i+v_0\not\in\Gamma_F\}.$$ In addition, he proves that $\Lambda_{F}\setminus\Gamma_{F}=\emptyset$ if and only if $\mathcal{C}_F$ is analytically equivalent to $\mathcal{C}_G$ with $G=y^{v_0}-x^{v_1}$ with $GCD(v_0,v_1)=1$. In this case we put $\lambda=\infty$.

The exponent $\lambda$ is called the {\bf Zariski invariant} of $\mathcal{C}_F$. Notice that if $\lambda\neq\infty$ then $v_1=\beta_1<\lambda \leq\beta_2$ and $\lambda+v_0$ is a minimal generator for $\Lambda_F$.

In what follows we present an alternative way to obtain the Zariski invariant of a plane branch using dicritical foliations in the first triple point of the dual graph $G(\mathcal{C}_F)$ as described in Corollary \ref{CorSep}.

\begin{lemma}\label{zinv}
If $\lambda$ is the Zariski invariant of a plane branch $\mathcal{C}_F$ with semigroup $\Gamma_{F}=\langle v_0,\ldots ,v_g\rangle$, then there exist $H_1, H_2\in \mathbb{C}\{x\}[y]$ with $deg_yH_l<n_1=\frac{v_0}{e_1}$, $H_1$ a unit and $v_0+v_1< I(F,H_2)$ such that $\nu(H_1\omega_{01}+dH_2)=\lambda+v_0$.
\end{lemma}
\Dem 	If $\lambda=\infty$ then $\Lambda_F=\Gamma_F\setminus\{0\}$ with $\Gamma=\langle v_0,v_1\rangle$, that is, $e_1=GCD(v_0,v_1)=1$ and in this case $v_0+v_1<\nu(\omega_{01})\in\Gamma_{F}\setminus\{0\}$. So, by  (\ref{economy}), there exists $G_1\in\langle x,y\rangle\cap\mathbb{C}\{x\}[y]$ with $deg_yG_1<v_0$ such that $\nu(\omega_{01})=\nu(dG_1)<\nu(\omega_{01}-dG_1)\in\Gamma_{F}$.

In the same way, we obtain $G_2,\ldots ,G_s\in\langle x,y\rangle\cap\mathbb{C}\{x\}[y]$  satisfying $\mu_F\leq\nu(\omega_{01}-\sum_{i=1}^{s}dG_i)\in\Gamma_F$ consequently, as $(\mathcal{O}:\overline{\mathcal{O}})=\langle t^{\mu_F}\rangle$,  there exists $G\in\langle x,y\rangle\cap\mathbb{C}\{x\}[y]$  such that $\varphi^*(dG)=\varphi^*(\omega_{01}-\sum_{i=1}^{s}dG_i)$, that is, $\omega=\omega_{01}-d(G+\sum_{i=1}^{s}G_i)$ satisfies $\varphi^*(\omega)=0$ or equivalently $\frac{\omega\wedge dF}{dx\wedge dy}\in\langle F\rangle$. So, the result is obtained taking $H_1=1$ and $H_2=-(G+\sum_{i=1}^{s}G_i)$

If $\infty\neq\lambda+v_0=\nu(\eta)$, then by Lemma \ref{azevedo}, we can express $\eta=A_1\omega_{01}+dA_2$ with $A_1, A_2\in\mathbb{C}\{x,y\}$. Considering a $2$-semiroot $F_2\in\mathbb{C}\{x\}[y]$ (recall that $F_2=F$ if $e_1=1$) we write $A_i=B_iF_2+H_i$, $deg_yH_i< deg_yF_2=\frac{v_0}{e_1}$ and
$$\eta=H_1\omega_{01}+dH_2+B_2dF_2+F_2\cdot (dB_2+B_1\omega_{01}).$$
As $\lambda\leq\beta_2$, it follows by (\ref{relations0}) that $\lambda+v_0< v_2\leq\nu(B_2dF_2+F_2\cdot (dB_2+B_1\omega_{01}))$. So, $\lambda+v_0=\nu(\zeta)$ where $\zeta:=H_1\omega_{01}+dH_2$ and, without loss of generality $H_2$ can be considered a non unit, that is, $H_2\in\langle x,y\rangle$.

We must have $I(F,H_2)> I(F,H_1)+v_0+v_1$. Indeed, if $I(F,H_2)\leq I(F,H_1)+v_0+v_1<\nu(H_1\omega_{01})$, then $\lambda+v_0=\nu(\zeta)=\nu(dH_2)\in\Gamma_F$ that is a contradiction, because $\lambda+v_0\in\Lambda_{F}\setminus\Gamma_{F}$.

Notice that $I(F,H_2)> v_0+v_1$ implies that $H_2 \in \langle y^2\rangle + \langle x,y\rangle^3$. In particular, $(H_2)_x \in \langle x, y\rangle^2$ and $(H_2)_y \in \langle x^2, y\rangle$.

Since $\zeta = (v_1yH_1+(H_2)_x)dx+(-v_0xH_1+(H_2)_y)dy$, if $H_1\in\langle x,y\rangle$ then $\zeta$ is expressed as $Mdx+Ndy$ with $M\in\langle x,y\rangle^2$ and $N\in\langle x^2,y\rangle$. In this way, by Remark \ref{elimination}, the exponent $\nu(\zeta)-v_0=\lambda$ could not be the Zariski invariant of $\mathcal{C}_F$. So, $H_1$ is a unit and we get the result considering $\zeta=H_1\omega_{01}+dH_2$.
\cqd

As $\lambda+v_0$ is a minimal generator for $\Lambda_{F}$, any $\omega\in\Omega^1$ such that $\nu(\omega)=\lambda+v_0$ defines a foliation. In this way, by the above lemma and Corollary \ref{CorSep}, we can compute the Zariski invariant for a plane branch $\mathcal{C}_F$ considering dicritical foliations in the first triple point of $G(\mathcal{C}_F)$ defined by $H_1\omega_{01}+dH_2$ with $H_1, H_2\in \mathbb{C}\{x\}[y]$ with $deg_yH_l<n_1=\frac{v_0}{e_1}$, $H_1$ a unit and $v_0+v_1< I(F,H_2)$, or equivalently, dicritical foliations defined by $w_{01}+H_1^{-1}dH_2=w_{01}+P_1dx+P_2dy$ with $\nu(P_1dx+P_2dy)=\nu(H_1^{-1}dH_2)=\nu(dH_2)=I(F,H_2)>v_0+v_1$.

Given $P_1dx+P_2dy\in\Omega^1$ and considering a $2$-semiroot $F_2\in\mathbb{C}\{x\}[y]$ of $F$ we may write
$$P_2=G_2\cdot (F_2)_y+Q_2\ \ \mbox{and}\ \ P_1-G_2\cdot (F_2)_x=G_1F_2+Q_1$$
with $G_1, G_2\in\mathbb{C}\{x,y\}$, $Q_1, Q_2\in\mathbb{C}\{x\}[y]$ with $deg_yQ_1<deg_yF_2=\frac{v_0}{e_1}$ and $deg_yQ_2<deg_y(F_2)_y=\frac{v_0}{e_1}-1$. In this way, we get
$$w_{01}+P_1dx+P_2dy=w_{01}+Q_1dx+Q_2dy+G_2dF_2+F_2G_1dx.$$
If $\nu(w_{01}+P_1dx+P_2dy)=\lambda+v_0$, then $\nu(Q_1dx+Q_2dy)>v_0+v_1$ and $\nu(w_{01}+Q_1dx+Q_2dy)=\lambda+v_0$, because $\nu(G_2dF_2+F_2G_1dx)\geq v_2>\beta_2+v_0\geq\lambda+v_0$.

Notice that the condition $\nu(Q_1dx+Q_2dy)>v_0+v_1$ is equivalent to consider $Q_1\in\langle x,y\rangle^2$, $mult_xQ_1(x,0)>\frac{v_1}{v_0}$ and $Q_2\in\langle x^2,y\rangle$. In this way, by the previous lemma and the above comments, $\lambda+v_0=\nu(\omega)$ with $\omega$ belonging to the set
$$\mathcal{D}_{1}=\left \{
\begin{array}{l}
	\omega_{01}+Q_1dx+Q_2dy;\ Q_1\in\langle x,y\rangle^2,\ Q_2\in\langle x^2,y\rangle\ \mbox{with}\\ deg_yQ_1<\frac{v_0}{e_1},\ deg_yQ_2<\frac{v_0}{e_1}-1\ \mbox{and}\ mult_xQ_1(x,0)>\frac{v_1}{v_0}
\end{array}\right \}.$$

Moreover, we have the following result:

\begin{proposition}\label{oziel-generalized}
For any plane branch $\mathcal{C}_F$ as considered in this subsection
$$\lambda=\max \{\nu(\omega)-v_0;\ \omega\in\mathcal{D}_1 \}=\max\left \{I\left(F,\frac{\omega\wedge dF}{dx\wedge dy}\right )-(\mu_F+v_0-1);\ \omega\in\mathcal{D}_1 \right \}.$$ In particular, $\lambda$ is determined considering foliations $\mathcal{F}_{\omega}$ with $\omega\in\mathcal{D}_1$.
\end{proposition}
\Dem By the above comments, there exists
$\omega=\omega_{01}+Q_1dx+Q_2dy\in\mathcal{D}_1$ such that $\nu(\omega)=\lambda+v_0$.

The case $\lambda=\infty$ is immediate. So, let us consider $\lambda\neq\infty$.

Suppose that there exists $\eta=Adx+Bdy\in\mathcal{D}_1$ such that $\nu(\eta)>\lambda+v_0$.

If $\nu(\eta)\in\Lambda_{F}\setminus\Gamma_{F}$, by  Proposition 1.3.13 in \cite{handbook}, we have that $A\in\langle x,y\rangle^2$ and $B\in\langle x^2,y\rangle$, but this contradicts the fact that $\eta\in\mathcal{D}_1$.

If $\nu(\eta)\in\Gamma_{F}$, by Lemma 1.3.12 in \cite{handbook}, there exists $A_1dx+B_1dy\in\Omega^1$ with $A_1\in\langle x,y\rangle^2$ and $B_1\in\langle x^2,y\rangle$ such that $\nu(\eta)=\nu(A_1dx+B_1dy)<\nu(\eta-A_1dx-B_1dy)$. Proceeding in this way we obtain $P_1\in\langle x,y\rangle^2$ and $P_2\in\langle x^2,y\rangle$
such that $\eta-P_1dx-P_2dy\in\mathcal{D}_1$ satisfies
$$\nu(\eta-P_1dx-P_2dy)\in\Lambda_{F}\setminus\Gamma_{F}\ \ \ \ \mbox{or}\ \ \ \ \nu(\eta-P_1dx-P_2dy)=\infty.$$
As before $\nu(\eta-P_1dx-P_2dy)\in\Lambda_{F}\setminus\Gamma_{F}$ contradicts the fact that $\eta-P_1dx-P_2dy\in\mathcal{D}_1$. On the other hand, if $\nu(\eta-P_1dx-P_2dy)=\infty$, as $\eta\in \mathcal{D}_1$ there exist $Q_1\in\langle x,y\rangle^2$ and $Q_2\in\langle x^2,y\rangle$ with $mult_x Q_1(x,0)>\frac{v_1}{v_0}$ such that $\eta=\omega_{01}+Q_1dx+Q_2dy$, then $\omega_{01}=(P_1-Q_1)dx+(P_2-Q_2)dy$ with $P_1-Q_1\in\langle x,y\rangle^2$ and $P_2-Q_2\in\langle x^2,y\rangle$, that is an absurd.

Hence, $\lambda+v_0=\max\{\nu(\omega);\ \omega\in\mathcal{D}_1\}$ and, by Remark \ref{jac}, $\nu(\omega)+\mu_F-1=\nu\left ( \frac{\omega\wedge dF}{dx\wedge dy}\right )=I\left ( F, \frac{\omega\wedge dF}{dx\wedge dy}\right )$ that concludes the proof.
\cqd

For plane branches with semigroup $\langle v_0,v_1\rangle$, the previous result was obtained by G\'omez-Mart\'inez in \cite{oziel}, where foliations defined by an element in $\mathcal{D}_1$ are called dicritical cuspidal foliations.

\section{Technical Lemmas}\label{technical}

In what follows, given $S(t)=\sum_{i\geq i_0}a_it^i\in\mathbb{C}\{t\}$ we denote ${\rm Coeff}(S(t),t^k):=a_k$, that is, the coefficient of $t^k$ in $S(t)$. In particular, if $H\in\mathbb{C}\{x,y\}$ and $\psi(t)=(t^{i_0},\sum_{i\geq i_1}d_it^i)$ then ${\rm Coeff}(\psi^*(H),t^k)$ depends polynomially on the coefficients $d_{i_1},\ldots ,d_{i_k}$ for some $i_k \geq i_1$. In this case, we write ${\rm Coeff}(\psi^*(H)),t^k)=p(d_{i_1},\ldots ,d_{i_k})\in\mathbb{C}[d_{i_1},\ldots ,d_{i_k}]$. If $i_k < i_1$ then we assume that $\mathbb{C}[d_{i_1},\ldots ,d_{i_k}]$ is $\mathbb{C}$.

Let $\mathcal{C}_F$ be a plane branch with semigroup $\Gamma_F=\langle v_0,v_1, \ldots, v_g \rangle$ and a canonical system of semiroots $\{F_0,F_1,\ldots ,F_g,F_{g+1}=F\}$ as in Proposition \ref{Prop-MinimalPolyn} with $$\varphi_{i}(t)=\left(t^{\frac{\beta_0}{e_{i-1}}},\sum_{\beta_1 \leq j <
	\beta_i}c_jt^{\frac{j}{e_{i-1}}}\right) \hspace{0.5cm} \mbox{and} \hspace{0.5cm} \varphi (t)=\left(t^{\beta_0},\sum_{j \geq \beta_1}c_jt^j\right)
$$ parameterizations of $\mathcal{C}_{F_i}$ $i=1,\ldots , g$ and $\mathcal{C}_F=\mathcal{C}_{F_{g+1}}$, respectively.

\begin{lemma}\label{coef-control} For each $k \geq v_i$ and  $i=1,\ldots, g$, we have ${\rm Coeff}(\varphi^*(F_i),t^k)\in\mathbb{C}[c_{\beta_1},\ldots ,c_{k-v_i+\beta_i}]$ of degree $1$ in $c_{k-v_i+\beta_i}$.
\end{lemma}
\Dem  For $1 \leq i \leq g$ we take $F_i$ as in (\ref{Fi}), so
$$\varphi^*(F_i)
	=\prod_{\alpha\in U_{m_i}}\left (\sum_{j\geq \beta_1
	}c_jt^j-\sum_{\beta_1\leq j<\beta_i}c_j\alpha^{\frac{j}{e_{i-1}}}
	t^j\right )=\prod_{\alpha\in U_{m_i}}\left (\sum_{\beta_1\leq
		j<\beta_i}c_j(1-\alpha^{\frac{j}{e_{i-1}}}) t^j+\sum_{j\geq \beta_i}c_jt^j\right ).$$

Denoting $G_s=\{\alpha \in \mathbb{C};\ \alpha^{\frac{e_s}{e_{i-1}}}=1\}$ for $0 \leq s < i$ we get
$$\{1\}=G_{i-1}\subset G_{i-2}\subset \cdots \subset G_1\subset G_0=\{\alpha  \in \mathbb{C};\ \alpha^{\frac{e_0}{e_{i-1}}}=1\}=U_{m_i}.$$

Notice that for each $j$ satisfying $\beta_0 \leq j < \beta_{s+1}$ and $\alpha \in G_s \setminus G_{s+1}$ we have $\alpha^{\frac{j}{e_{i-1}}}=1$, since $e_s$ divides $j$. Thus
$$\sum_{\beta_1\leq j<\beta_i}c_j(1-\alpha^{\frac{j}{e_{i-1}}}) t^j+\sum_{j\geq \beta_i}c_jt^j=\left \{\begin{array}{ll}
	c_{\beta_{s+1}}(1-\alpha^{\frac{\beta_{s+1}}{e_{i-1}}})t^{\beta_{s+1}}+ h.o.t. & \mbox{if}\ \alpha\in G_{s}\setminus G_{s+1} \vspace{0.2cm}\\
	c_{\beta_i}t^{\beta_i}+ h.o.t. & \mbox{if}\ \alpha\in G_{i-1}=\{1\},
\end{array}
\right .$$ with $\alpha^{\frac{\beta_{s+1}}{e_{i-1}}}\neq 1$ if
$\alpha\in G_s\setminus G_{s+1}$ (see \cite{hefez}, Lemma 6.8).

It follows that
\begin{equation}\label{I(F,Fi)}
	\begin{array}{rcl}
		\varphi^*(F_i) & = & \displaystyle\left( \prod_{l=0}^{i-2}\ \prod_{\alpha\in G_{l}\setminus G_{l+1}}\left ( (1-\alpha^{\frac{\beta_{l+1}}{e_{i-1}}})c_{\beta_{l+1}}t^{\beta_{l+1}}+ h.o.t. \right)\right )\cdot\prod_{\alpha \in G_{i-1}}(c_{\beta_i}t^{\beta_i}+ h.o.t.)\vspace{0.2cm}\\
		&=& \displaystyle \prod_{l=0}^{i-2}\ \left(b_l t^{\frac{e_l-e_{l+1}}{e_{i-1}}\beta_{l+1}} +h.o.t.\right)\cdot \left(c_{\beta_i}t^{\beta_i}+ h.o.t.\right)=\left(\prod_{l=0}^{i-2}b_l\right)c_{\beta_i}t^{v_i}+h.o.t.,
	\end{array}
\end{equation}
where $b_l:=\prod_{\alpha\in G_{l}\setminus G_{l+1}}(1-\alpha^{\frac{\beta_{l+1}}{e_{i-1}}})c_{\beta_{l+1}} \neq 0$ for $l=0,\ldots, i-2$ since $c_{\beta_i} \neq 0$, it is easy to verify that $\sharp (G_s\setminus G_{s+1})=\frac{e_s-e_{s+1}}{e_{i-1}}$ for $0\leq s\leq i-1$, and the last equality follows by (\ref{relations}) since $\sum_{l=0}^{i-2}\frac{e_l-e_{l+1}}{e_{i-1}}\beta_{l+1}+\beta_i=v_i$. In fact, this is a proof that $I(F,F_i)=v_i$.

In addition, by the previous analysis, we obtain
${\rm Coeff}(\varphi^*(F_i),t^k)\in\mathbb{C}[c_{\beta_1},\ldots ,c_{\rho}]$ with
$k=\frac{1}{e_{i-1}}\sum_{l=0}^{i-2}(e_l-e_{l+1})\beta_{l+1}+\rho$,
that is, $\rho=k-v_i+\beta_i$ and ${\rm Coeff}(\varphi^*(F_i),t^k)$ is a polynomial of degree one in $c_{k-v_i+\beta_i}$. More explicitly, we get \begin{equation}\label{coef-Fi}
	{\rm Coeff}(\varphi^*(F_i), t^k)=P_{ik}(c_{\beta_1},\ldots
	,c_{k-v_i+\beta_i-1})+\left ( \prod_{l=0}^{i-2}b_l\right
	)c_{k-v_i+\beta_i} \end{equation} where
$P_{ik}(c_{\beta_1},\ldots
,c_{k-v_i+\beta_i-1})$ is a homogeneous polynomial.
\cqd

In what follows we consider the family $\mathcal{C}_{F_a}$ of plane branches topologically equivalent to $\mathcal{C}_{F}$ parameterized as (\ref{family}), that is
$$\varphi_a(t)=\left (t^{\beta_0},\sum_{\beta_1\leq l<\beta_{g}}c_lt^l+\sum_{l\geq\beta_{g}}a_lt^l\right ),$$
where $c_l$ is the coefficient of $t^l$ in $\varphi(t)$, then it is constant  for $\beta_1\leq l<\beta_{g}$ and $a_l$ is a complex parameter for $l\geq\beta_{g}$ with $a_{\beta_{g}}$ nonvanishing.

It is immediate that $\Gamma_{F_a}=\Gamma_{F}=\langle
v_0,v_1,\ldots ,v_g\rangle$,
$\{F_0,F_1,\ldots ,F_{g},F_a\}$ and $\{F_0,F_1,\ldots ,F_{g},F\}$ are the canonical system of semiroots
of $F_a$ and $F$, respectively. In particular, $I\left ( F,\prod_{i=0}^{g} F_i^{\gamma_i} \right )=I\left ( F_a,\prod_{i=0}^{g} F_i^{\gamma_i}  \right )$ for any nonnegative integers $\gamma_i$.

\begin{remark}\label{coefFi}
	As $c_l$ is constant for $\beta_1\leq l<\beta_g$, by Lemma \ref{coef-control} or more precisely by (\ref{coef-Fi}), we obtain $${\rm Coeff}(\varphi^*_a(F_i),t^k)=p_{ik}(a_{\beta_g},\ldots ,a_{k-v_i+\beta_i-1})+\delta_{ik}\cdot a_{k-v_i+\beta_i}$$ with $\delta_{ik}\in\mathbb{C}^*$ for any $k\geq v_i$ and $1\leq i\leq g$.
	In particular, we get
	$${\rm Coeff}(\varphi^*_a(dF_i),t^k)=(k+1)\cdot (p_{i,k+1}(a_{\beta_g},\ldots ,a_{k-v_i+\beta_i})+ \delta_{i,k+1}\cdot a_{k-v_i+\beta_i+1}).$$	
	
	Moreover, for $k<v_i+\beta_g-\beta_i$ we have ${\rm Coeff}(\varphi^*_a(F_i),t^k)$ is constant, that is, it does not depend on the parameters $a_l$ for $l \geq \beta_g$,  and by (\ref{I(F,Fi)})
	$$
		{\rm Coeff}(\varphi^*_a(F_i),t^{v_i})=\left \{
		\begin{array}{ll}
			\delta_{i} & \mbox{if}\ i<g \\
			\delta_{g} \cdot a_{\beta_g} & \mbox{if}\ i=g,\end{array}\right .
$$ for some $\delta_{i}, \delta_{g}\in\mathbb{C}^*$.
\end{remark}

The next lemmas are variations of Lemma \ref{coef-control}.

\begin{lemma}\label{coef-control-bg} Given $H=\prod_{i=0}^{j}F_i^{\gamma_i}$ with $\gamma_j\neq 0$ for some $0 \leq j \leq g$, we get ${\rm Coeff}(\varphi^*_a(H),t^k)\in\mathbb{C}[a_{\beta_g},\ldots,a_{k-I(F,H)+\beta_j}]$ for $k\geq I(F,H)$.
	In particular,
	$${\rm Coeff}(\varphi^*_a(H),t^{I(F,H)})=\left \{
	\begin{array}{ll}
		\delta_{jH} & \mbox{if}\ j<g \\
		\delta_{gH}\cdot a_{\beta_g}^{\gamma_g} & \mbox{if}\ j=g,\end{array}\right .$$
	for some $\delta_{jH},\delta_{gH}\in\mathbb{C}^*$.
\end{lemma}
\Dem Since $F_0=x$, if $j=0$ then $H=F_0^{\gamma_0}$ and $I(F,H)=I(F_a,H)=\gamma_0\cdot v_0$ with
$\varphi^*_a(H)=\varphi^*_a(F_0^{\gamma_0})=t^{\gamma_0v_0}$ and we get the result.

Notice that for any $\gamma_l>0$ with $l\neq 0$, we may rewrite
\begin{equation}
	\label{H-Semirrots}
	\varphi^*_a(H)=\varphi^*_a(F_0^{\gamma_0})\cdot \varphi^*_a(F_1^{\gamma_1})\cdot\ldots\cdot\underbrace{\varphi^*_a(F_l)\cdot \varphi^*_a(F_l^{\gamma_l-1})}\cdot\ldots\cdot \varphi^*_a(F_j^{\gamma_j}).
\end{equation}

In order to determine $\epsilon\in\mathbb{N}$ such that ${\rm Coeff}(\varphi_{a}^*(H)),t^k)\in\mathbb{C}[a_{\beta_g},\ldots ,a_{\epsilon}]$ it is sufficient to analyse for each $l$, the product of a term of order $\gamma$ in $\varphi_{a}^*(F_l)$ and all the initial terms of factors in (\ref{H-Semirrots}) such that $k=\gamma +\sum_{i=0}^j\gamma_iv_i - v_l$, that is, $\gamma=k-\sum_{i=0}^{j}\gamma_i v_i+v_l=k-I(F,H)+v_l$, since $I(F,H)=I(F_a,H)$.

In this way, we can determine $\epsilon$ considering
\begin{equation}\label{C}
\prod_{i=0\atop i\neq l}^{j}\left ({\rm Coeff}(\varphi^*_{a}(F_i),t^{v_i})\right )^{\gamma_i}\cdot \left ({\rm Coeff}(\varphi^*_{a}(F_l),t^{v_l})\right )^{\gamma_{l}-1}\cdot {\rm Coeff}(\varphi^*_{a}(F_l),t^{k-I(F,H)+v_l}).\end{equation}

By Remark \ref{coefFi}, the expression (\ref{C}) is polynomial in $\mathbb{C}[a_{\beta_g},\ldots ,a_{\gamma-v_l+\beta_l}]=\mathbb{C}[a_{\beta_g},\ldots ,a_{k-I(F,H)+\beta_l}]$ of degree one in $a_{k-I(F,H)+\beta_l}$. Thus
$$\epsilon=\max_{0\leq l\leq j}\{k-I(F,H)+\beta_l;\ \gamma_l\neq 0\}
=k-I(F,H)+\beta_j.$$

By (\ref{C}) and the previous remark, we have
$${\rm Coeff}(\varphi^*_{a}(H),t^k)=p_{kH}(a_{\beta_g},\ldots ,a_{k-I(F,H)+\beta_j-1})+\left \{
\begin{array}{ll}
	\delta_{kH}\cdot a_{k-I(F,H)+\beta_j} & \mbox{if}\ j<g \vspace{0.2cm}\\
	\delta_{kK}\cdot a_{\beta_g}^{\gamma_g-1}\cdot a_{k-I(F,H)+\beta_g} & \mbox{if}\ j=g
\end{array}\right .
$$
for some $\delta_{kH}\in\mathbb{C}^*$.
In particular,
$${\rm Coeff}(\varphi^*_a(H),t^{I(F,H)})=\left \{
\begin{array}{ll}
	\delta_{jH} & \mbox{if}\ j<g \\
	\delta_{gH} \cdot a_{\beta_g}^{\gamma_g} & \mbox{if}\ j=g,\end{array}\right .$$
for some $\delta_{jH},\delta_{gH}\in\mathbb{C}^*$.
\cqd

As an immediate consequence of Lemma \ref{coef-control-bg} we obtain the following result.

\begin{corollary}\label{coef-dH} If $H=\prod_{i=0}^{j}F_i^{\gamma_i}$ with $\gamma_j\neq 0$ for some $0 \leq j \leq g$, then for any $k> I(F,H)$ we get ${\rm Coeff}(\varphi^*_a(dH),t^k)\in\mathbb{C}[a_{\beta_g},\ldots ,a_{k-I(F,H)+\beta_j+1}]$.
\end{corollary}

Now we analyse  $\varphi^*_{a}(\omega_{ig})$ for a $1$-form
$\omega_{ig}=v_iF_idF_g - v_gF_g dF_i,$  considered in (\ref{wij}) with $0 \leq i < g$.
Notice that, by Remark \ref{coefFi}, we get  $ord_t\varphi^*_{a}(\omega_{ig})\geq v_i+v_g$.

\begin{lemma}\label{coef-control-wij} For $k \geq v_i+v_g$ we have
	$${\rm Coeff}(\varphi^*_a(\omega_{ig}),t^k)=q_{ik}(a_{\beta_g},\ldots
	,a_{k-v_i-v_g+\beta_g})+r_{ik}\cdot a_{k-v_i-v_g+\beta_g+1} ,$$ with $r_{ik} \in \mathbb{C}^*$.
\end{lemma}
\Dem The highest order of a term in $\varphi_a(t)$ that contributes
with the term of order $k$ in
$v_i\varphi_a^*(F_i)\varphi_a^*(dF_g)$ is determined by
considering the product of the initial term $\delta_{i}\cdot t^{v_i}$ of $\varphi_a^*(F_i)$ with the term of order $k-v_i$ in
$\varphi_a^*(dF_g)$ or the product of the initial term $v_g\cdot \delta_{g}\cdot a_{\beta_g}t^{v_g-1}$ of $\varphi_a^*(dF_g)$ with the term with order $k-v_g+1$ in
$\varphi_a^*(F_i)$.

By Remark \ref{coefFi}, ${\rm Coeff}(\varphi_{a}^*(F_i),t^{k-v_g+1})\in\mathbb{C}[a_{\beta_g},\ldots ,a_{k-v_i-v_g+\beta_i+1}]$ and ${\rm Coeff}(\varphi_a^*(dF_g),t^{k-v_i})=(k-v_i+1)\cdot (p_{g,k-v_i+1}(a_{\beta_g},\ldots ,a_{k-v_i-v_g+\beta_g})+ \delta_{g,k-v_i+1}\cdot a_{k-v_i-v_g+\beta_g+1})$. In this way,
$$
{\rm Coeff}(\varphi_a^*(F_i),t^{v_i})\cdot {\rm Coeff}(\varphi_a^*(dF_g),t^{k-v_i})\in\mathbb{C}[a_{\beta_g},\ldots ,a_{k-v_i-v_g+\beta_g+1}]$$
and ${\rm Coeff}(\varphi_a^*(dF_g),t^{v_g-1})\cdot {\rm Coeff}(\varphi_{a}^*(F_i),t^{k-v_g+1}) \in\mathbb{C}[a_{\beta_g},\ldots ,a_{k-v_i-v_g+\beta_i+1}].$

Moreover, with the notation of Remark \ref{coefFi} and denoting $e:=\delta_{i}\cdot\delta_{g,k-v_i+1} \neq 0$, we get
$${\rm Coeff}(v_i\varphi_a^*(F_idF_g),t^k)=Q_{1k}(a_{\beta_g},\ldots ,a_{k-v_i-v_g+\beta_g})+v_i\cdot (k-v_i+1)\cdot e\cdot a_{k-v_i-v_g+\beta_g+1}.$$

In a similar way, ${\rm Coeff}(v_g\varphi_a^*(F_gdF_i))=Q_{2k}(a_{\beta_g},\ldots ,a_{k-v_i-v_g+\beta_g})+v_g\cdot v_i\cdot e\cdot a_{k-v_i-v_g+\beta_g+1}.$ Hence,
$$	{\rm Coeff}(\varphi_{a}^*(\omega_{ig}),t^k)=q_{ik}(a_{\beta_g},\ldots
	,a_{k-v_i-v_g+\beta_g})+ r_{ik} \cdot a_{k-v_i-v_g+\beta_g+1},
$$ where
$q_{ik}=Q_{1k}-Q_{2k}\in\mathbb{C}[a_{\beta_g},\ldots
,a_{k-v_i-v_g+\beta_g}]$ and $r_{ik}:=e \cdot v_i\cdot (k-v_i-v_g+1)\neq 0.$\cqd

In the next result, we determine the coefficient of a term in $\varphi_a^*(H\cdot\omega_{ig})$, for any $H \in \mathbb{C}\{x\}[y]$ expressed according to Proposition \ref{Prop-DecompH}, that is,
$H=\sum_{\delta}e_{\delta}F_0^{\delta_0}F_1^{\delta_1}\cdots F_g^{\delta_g}$  with $e_{\delta}\in\mathbb{C}^*$ and such that
$I(F,H)=I(F,F_0^{\gamma_0}F_1^{\gamma_1}\cdots
F_g^{\gamma_g}),$ for some non-negative integers $\gamma_0, \gamma_1,\ldots, \gamma_g$.

\begin{lemma}\label{coef-control-Hwij}
	If $H=\sum_{\delta}e_{\delta}F_0^{\delta_0}F_1^{\delta_1}\cdots F_g^{\delta_g}$ is as in
	(\ref{semiroots}) with
	$I(F,H)=I(F,F_0^{\gamma_0}F_1^{\gamma_1}\cdots
	F_g^{\gamma_g}),$ then for $k \geq I(F,H)+v_i+v_g$  we get
	$${\rm Coeff}(\varphi_a^*(H\cdot\omega_{ig}),t^k)=p_{k}(a_{\beta_g},\ldots,a_{k_{ig}-1})+r_{k}\cdot
	a_{\beta_g}^{\gamma_g}\cdot a_{k_{ig}}$$
	with $k_{ig}:=k-I(F,H)-v_i-v_g+\beta_g+1$ and some $r_{k}\in\mathbb{C}^*$.
\end{lemma}
\Dem The coefficient of $t^k$ in $\varphi_a^*(H\cdot\omega_{ig})$ is obtained by the sum of products of a term of order $s_1$ in $\varphi_a^*(H)$ and a term of order $s_2$ of $\varphi_a^*(\omega_{ig})$ such $s_1+s_2=k$. In this way, to prove the lemma it is sufficient to analyse a such product for the maximum possible value for $s_1$ or $s_2$.

For each element $e_{\delta}F_0^{\delta_0}F_1^{\delta_1}\cdots F_g^{\delta_g}$ we set \begin{equation}
	\label{mdelta}
	m_{\delta}:=\max_{0\leq l\leq g}\{l;\ \delta_l\neq 0\}\ \ \mbox{and}\ \ I_{\delta}:=I\left (F,\prod_{l=0}^{g}F_l^{\delta_l}\right )=I\left (F_a,\prod_{l=0}^{g}F_l^{\delta_l}\right ).
\end{equation}

\vspace{0.2cm}

\noindent{\bf Case 1)} $s_1=I(F,H)=I(F_a,H)$ and $s_2=k-I(F,H)$, that is, the maximum possible value for $s_2$.

\vspace{0.1cm}
Notice that ${\rm Coeff}(\varphi_{a}^*(H),t^{I(F,H)})={\rm Coeff}(\prod_{l=0}^{g}\varphi_{a}^*(F_l^{\gamma_l}),t^{I(F,H)})$ and, by Lemma \ref{coef-control-bg}, it is constant if $m_{\gamma}<g$ and equal to $e_H\cdot a_{\beta_g}^{\gamma_g}$ with $e_H\in\mathbb{C}^*$ for $m_{\gamma}=g$. On the other hand, by previous lemma ${\rm Coeff}(\varphi_{a}^*(\omega_{ig}),t^{k-I(F,H)})=p(a_{\beta_g},\ldots
,a_{k_{ig}-1})+r\cdot a_{k_{ig}}$ where $k_{ig}:=k-I(F,H)-v_i-v_g+\beta_g+1$ and $r\in \mathbb{C}^*$. Hence,
${\rm Coeff}(\varphi_{a}^*(H),t^{I(F,H)})\cdot {\rm Coeff}(\varphi_{a}^*(\omega_{ig}),t^{k-I(F,H)})$ is expressed as \begin{equation}\label{caso1}
	p_k(a_{\beta_g},\ldots
	,a_{k_{ig}-1})+r_k\cdot a_{\beta_g}^{\gamma_g}\cdot a_{k_{ig}}
\end{equation} with $r_k:=r \cdot e_H \neq 0$.

\noindent{\bf Case 2)} $s_1=k-v_i-v_g$ and $s_2=v_i+v_g$ that is, the maximum possible value for $s_1$.

\vspace{0.1cm}
As $k\geq I(F,H)+v_i+v_g$, if $k=I(F,H)+v_i+v_g$ then we must have $s_1=I(F,H)$ and $s_2=k-I(F,H)$, that is, we are in the previous case. So, we can assume that $k>I(F,H)+v_i+v_g$.

Notice that ${\rm Coeff}(\varphi_{a}^*(H),t^{k-v_i-v_g})=\sum_{\delta}e_{\delta}{\rm Coeff}(\prod_{l=0}^{g}\varphi_{a}^*(F_l^{\delta_l}),t^{k-v_i-v_g})$. By Lemma \ref{coef-control-bg} and by (\ref{mdelta}), we get
$${\rm Coeff}\left (\prod_{l=0}^{g}\varphi_{a}^*(F_l^{\delta_l}),t^{k-v_i-v_g}\right )\in\mathbb{C}[a_{\beta_g},\ldots ,a_{k-v_i-v_g-I_{\delta}+\beta_{m_{\delta}}}].$$ Moreover, by Lemma \ref{coef-control-wij}, we have ${\rm Coeff}(\varphi_{a}^*(\omega_{ig}),t^{v_i+v_g})\in\mathbb{C}[a_{\beta_g},a_{\beta_g+1}].$ As  $I_{\delta}\geq I(F,H)$ and $k> I(F,H)+v_i+v_g$ we have that
$$
\max_{\delta}\{k-v_i-v_g-I_{\delta}+\beta_{m_{\delta}},\beta_{g}+1\}<k_{ig}.$$

Considering the above cases we conclude that ${\rm Coeff}(\varphi_a^*(H\cdot\omega_{ig}),t^k)$ is  given as (\ref{caso1}).  \cqd

\vspace{0.5cm}

{\sc Nuria Corral}. Departamento de Matemáticas, Estadística y Computación, Universidad de Cantabria.   Av. de los Castros s/n, 39005, Santander, Spain. Email address: {\it nuria.corral@unican.es}\\

{\sc Marcelo E. Hernandes}. Departamento de Matem\'atica, Universidade Estadual de Maring\'a. Av. Colombo, 5790, 87020-900, Maring\'a - PR, Brazil. Email address: {\it mehernandes@uem.br}\\

{\sc Maria Elenice R. Hernandes}. Departamento de Matem\'atica, Universidade Estadual de Maring\'a. Av. Colombo, 5790, 87020-900, Maring\'a - PR, Brazil. Email address: {\it merhernandes@uem.br}


\begin{thebibliography}{30}

\bibitem{Abh} {\sc Abhyankar, S. S. and Moh, T.}, {\it Newton-Puiseux expansion and generalized Tschirnhausen transformation.} J. Reine Angew. Math., 260, 47-83 (1973); 261, 29-54, (1973).
	
\bibitem{azevedo} {\sc Azevedo, A. C. P.},
	{\it The Jacobian ideal of a plane algebroid curve}.
	Ph. D. Thesis, Purdue University (1967).
	
	
\bibitem{bayer} {\sc Bayer, V. and Hefez, A.}, {\it Algebroid plane curves whose Milnor and Tjurina numbers differ by one or two}. Bull. Braz. Math. Soc., 32(1), 63-81, (2001).
	
\bibitem{Cam-LN-S}  {\sc Camacho, C.;  Lins Neto, A. and  Sad, P.}: \textit{Topological invariants and equidesingularisation for holomorphic
vector fields.} J. Differential Geometry, \textbf{20}, 1 , 143--174 (1984).

\bibitem{cano} {\sc Cano, F., Cerveau, D. and D\'eserti, J.}, {\it Th\'eorie \'El\'ementaire des Feuilletages Holomorphes Singuliers}. Collection \'Echelles, Belin, (2013).
	
\bibitem{Can-C-2006} {\sc Cano, F., Corral, N.}, {\em Dicritical Logarithmic Foliations.} Publ. Mat. 50, 87–102, (2006).

\bibitem{Can-C-2011}  {\sc Cano, F., Corral, N.}, {\em Absolutely Dicritical Foliations.}  International Mathematics Research Notices, No. 8, 1926–1934, (2011).

\bibitem{Can-C-M-2019} {\sc Cano, F.;  Corral, N. and  Mol, R.}, \textit{Local polar invariants for plane singular foliations.\/}  Expo. Math.  \textbf{37},  no. 2, 145–164 (2019).

\bibitem{nuria} {\sc Cano, F.; Corral, N. and Senovilla-Sanz, D.}, {\it Analytic semiroots for plane branches and singular foliations}.  Bull. Braz. Math. Soc. (N.S.) 54, no. 2, Paper No. 27, 49 pp. (2023).

\bibitem{Cor-2024} {\sc Corral, N.}, \textit{Jacobian and polar curves of singular foliations. } To Appear in ``Handbook of Geometry and Topology of Singularities V: Foliations''.
	
	\bibitem{oziel} {\sc G\'omez-Mart\'inez, O.}, {\it Zariski invariant for non-isolated separatrices through Jacobian curves of pseudo-cuspidal dicritical foliations}. Journal of Singularities, 23, 236-270, (2021).
	
	\bibitem{hefez} {\sc Hefez, A.}, {\it Irreducible Plane Curve Singularities}. Real and Complex Singularities. Eds D. Mond and M. J. Saia. Lecture Notes in Pure and Appl. Math., 232, Marcel Dekker, N.Y., (2003).
	
	\bibitem{coloquio} {\sc Hefez, A. and Hernandes,
		M. E.}, {\it Computational Methods in the Local Theory of Curves}.  Mathematical Publications. 23o. Colóquio Brasileiro de Matemática. Instituto de Matemática Pura e Aplicada (IMPA), Rio de Janeiro, (2001).
	
	\bibitem{basestandard} {\sc Hefez, A. and Hernandes,
		M. E.}, {\it Standard bases for local rings of branches and their module of differentials}. J. Symb. Comp. 42, 178-191, (2007).
	
	\bibitem{HH} {\sc Hefez, A. and Hernandes, M. E.}, {\it The analytic classification of plane branches}. Bull. London Math. Soc., 43(2), 289-298, (2011).
	
	\bibitem{handbook} {\sc Hefez, A. and Hernandes, M. E.}, {\it The analytic classification of irreducible plane curve singularities}. Handbook of Geometry and Topology of Singularities II. Eds Cisneros-Molina, J. L., L\^e, D. T., Seade, J., Springer, 1-65, (2021).
	
	\bibitem{loray} {\sc Loray, F.}, {\it R\'{e}duction formelle des singularit\'{e}s cuspidales de champs de vecteurs analytiques}. Journal of Differential Equations 158, 152-173, (1999).
	
	\bibitem{Mat-S} {\sc Mattei, J.-F. and Salem, E.}, {\it Modules formels locaux de feuilletages holomorphes}. ArXiv:math/0402256.
	
\bibitem{Merle} {\sc Merle, M.}, {\it Invariants polaires des courbes planes}. Inventiones Mathematicae, 41, 103-112, (1977).	
	
	\bibitem{pol} {\sc Pol, D.}, {\it On the values of logarithmic residues along
		curves}.  Ann. L'Institut Fourier, 68(2), 725-766, (2018).
	
	\bibitem{Popescu} {\sc Popescu-Pampu, P.}, {\it Approximate roots.}  Valuation Theory and its Applications 33, AMS, 285-321, (2003).
	
	\bibitem{wall} {\sc Wall, C. T. C.}, {\it Singular Points of Plane Curves}. Cambridge University Press, (2004).
	
	\bibitem{torsion} {\sc Zariski, O.}, {\it Characterization of plane algebroid curves whose module of differentials has maximum torsion}. Proc. Nat. Acad. Sc. U.S.A. 56, 781-786, (1966).
	
	\bibitem{Zariski} {\sc Zariski, O.}, {\it Le Probl\`{e}me des Modules pour le Branches Planes.} Cours donn\'{e} au Centre de Math\'{e}matiques de L'\'{E}cole Polytechnique. English translation by Ben Lichtin: ``The Moduli Problem for Plane Branches'', University Lecture Series, 39, AMS, (2006).
	
\end{thebibliography}
\end{document}